\documentclass[12pt,a4paper]{article}
\usepackage[intlimits]{amsmath}
\usepackage{amsfonts,amssymb,amscd,amsthm}
\usepackage{cite}

\theoremstyle{plain}
\newtheorem{theorem}{Theorem}
\newtheorem{lemma}[theorem]{Lemma}
\newtheorem{proposition}[theorem]{Proposition}
\theoremstyle{definition}
\newtheorem{definition}[theorem]{Definition}

\theoremstyle{remark}
\newtheorem{remark}[theorem]{Remark}

\def\RR{\mathbb{R}}

\def\bcdot{\,\boldsymbol\cdot\,}
\def\lra{\longrightarrow}

\def\a{\alpha}
\def\be{\beta}
\def\ga{\gamma}
\def\dt{\delta}
\def\la{\lambda}
\def\e{\varepsilon}
\def\si{\sigma}
\def\ph{\varphi}
\def\om{\omega}

\def\Om{\Omega}
\def\Si{\Sigma}

\DeclareMathOperator\im{Im}
\DeclareMathOperator*\slim{s-lim}
\DeclareMathOperator*\esssup{ess\,sup}
\DeclareMathOperator\supp{supp}
\DeclareMathOperator\id{id}

\newcommand{\op}[1]{\operatorname{#1}}
\newcommand{\Vect}{{\rm Vect}}

\def\cA{{\mathcal{A}}}
\def\cB{{\mathcal{B}}}
\def\cC{{\mathcal{C}}}
\def\cF{{\mathcal{F}}}
\def\cI{{\mathcal{I}}}
\def\cJ{{\mathcal{J}}}
\def\cK{{\mathcal{K}}}
\def\cL{{\mathcal{L}}}
\def\cM{{\mathcal{M}}}
\def\cN{{\mathcal{N}}}
\def\cV{{\mathcal{V}}}
\def\gH{{\mathfrak{H}}}
\def\gZ{{\mathfrak{Z}}}

\def\ov{\overline}
\def\ovs{\overset}
\def\wt{\widetilde}
\def\wh{\widehat}
\def\pa{\partial}

\newcommand\pd[2]{\frac{\pa #1}{\pa #2}}

\let\rom\textup

\newcommand{\BL}{\biggl}
\newcommand{\BR}{\biggr}
\newcommand{\Bl}{\Bigl}
\newcommand{\Br}{\Bigr}
\newcommand{\bl}{\bigl}
\newcommand{\br}{\bigr}

\newcommand{\abs}[1]{\lvert#1\rvert}
\def\ka{\varkappa}
\newcommand{\norm}[1]{\left\| #1 \right\|}
\newcommand{\bnorm}[1]{ \pmb\|#1\pmb\|}

\numberwithin{equation}{section}

\let\thebibbliography\thebibliography

\def\thebibliography#1{\thebibbliography{#1}\addcontentsline{toc}{section}{Bibliography}}

\title{\bf Pseudodifferential Operators\\ on Stratified Manifolds}

\author{\bf V.~E.~Nazaikinskii, A.~Yu.~Savin, and B.~Yu.~Sternin}

\date{}

\begin{document}
\maketitle

\tableofcontents

\section*{Introduction}
\addcontentsline{toc}{section}{Introduction}

Pseudodifferential operators ($\Psi$DO) are an important tool in elliptic
theory on manifolds with singularities. When studying ``topological''
issues of this theory, including the homotopy classification and index
theory of elliptic operators, the main question of interest in $\Psi$DO
theory is the structure of the algebra of principal symbols (which
determine the Fredholm property of an operator) and the relationship
between this algebra and the algebra of operators themselves, in
particular, the compatibility conditions for the principal symbol
components corresponding to various strata of the manifold and the
existence of a quantization procedure (i.e., the construction of an
operator from its principal symbol). More subtle issues related to the
so-called complete symbols of $\Psi$DO and specific analytic formulas
describing $\Psi$DO do not play any significant role in this field.

In this connection, it is of interest to devise a $\Psi$DO theory in which
the principal symbol plays the main role and lower-order terms are touched
only when this is absolutely necessary. Such a construction is given in the
present paper.

We consider a specific class of manifolds with singularities, namely,
stratified manifolds, and describe a class of $\Psi$DO related to
differential operators with degeneration of first-order with respect to the
distance to the strata on such manifolds. (For manifolds with isolated
singularities, this is a Fuchs type degeneration.)

We restrict ourselves to the case of zero-order $\Psi$DO in $L^2$ spaces,
which is of main interest here. (Operators of nonzero order can be reduced
to zero-order operators by the order reduction procedure.) The
definition of $\Psi$DO adopted in our approach uses induction over the
number of strata in the manifold. The induction process inevitably involves
$\Psi$DO with parameters; hence we state all definitions for $\Psi$DO with
parameters from the very beginning. This simplifies the argument greatly
and does not lead to any serious complications in the statements
themselves.

Note that $\Psi$DO on manifolds with singularities were considered by
numerous authors. For example, $\Psi$DO with cone and edge degeneration in
so-called edge Sobolev spaces were considered in \cite{Schu1,EgSc1}; the
methods of these papers were developed for more complicated singularities
in~\cite{CMSch1}; a close class of $\Psi$DO on manifolds with corners was
studied in~\cite{Mel3}. Groupoids were used in \cite{Nis99} to analyze
$\Psi$DO on manifolds with singularities, and the paper~\cite{PlSe6} uses
localization methods in $C^*$-algebras for the same purpose.

In this paper, we first develop some auxiliary results pertaining to the
localization principle for abstract local operators (concerning this
principle, see also~\cite{GoKru2,AnLe1,PlSe6,Vas3,Vas1,Sim1}). Then we
consider $\Psi$DO.

Some results were partly announced without proof in our
paper~\cite{NSaSSc100} written jointly with B.-W.~Schulze.

\section{Generalized elliptic operators}

Let $H$ be a Hilbert space. By $\cB H$ and $\cK H$ we denote the algebra of
bounded linear operators in $H$ and the ideal of compact operators,
respectively. Suppose that $H$ is equipped with the structure of a unital
$*$-module over the $C^*$-algebra $C(X)$ of continuous functions on a
compact set $X$. An operator $A\in\cB H$ is called a \textit{generalized
elliptic operator} in the sense of Atiyah~\cite{Ati4} if it is Fredholm and
\textit{local}, i.e., compactly commutes with the action of $C(X)$. We need
a natural generalization of these notions to the case of operator
families.

\subsection{Local operators with parameter}

Let $Y$ be a Hausdorff locally compact topological space, in general,
noncompact (which will play the role of a parameter space). The
$C^*$-algebra $C(Y,\cB H)$ of bounded continuous operator families
$A:Y\lra\cB H$ with the $\sup$-norm
\begin{equation*}
    \bnorm{A}=\sup_{y\in Y}\norm{A(y)}
\end{equation*}
contains the closed $*$-ideal
\begin{equation*}
    \cJ=C_0(Y,\cK H)
\end{equation*}
of compact-valued families tending to zero at infinity in norm. By $\cA$ we
denote the $C^*$-subalgebra of $C(Y,\cB H)$ consisting of families $A(y)$
such that
\begin{equation}\label{commut}
    [A(y),\ph]\in\cJ\quad\text{for any function $\ph\in C(X)$.}
\end{equation}
The elements of $\cA$ are called \textit{local operators with a parameter}.

\subsection{Ellipticity and the Fredholm property}

The Fredholm property can naturally be generalized to the case of operator
families as follows. A family $A\in C(Y,\cB H)$ is said to be
\textit{Fredholm with parameter $y$} if the operator $A(y)$ is Fredholm for
all $y$ and invertible for large $y$ (i.e., outside some subset in $Y$) and
the inverse $A^{-1}(y)$ is bounded uniformly with respect to $y$ for large
$y$.

Now we can give the notion of generalized ellipticity in the sense of Atiyah.
An operator family $A\in C(Y,\cB H)$ is called a \textit{generalized elliptic
operator with a parameter} if it is local (i.e., $A\in\cA$) and Fredholm
with parameter in the sense described above.

Let
\begin{equation*}
 \si:\cA\lra\cA\slash\cJ
\end{equation*}
be the natural projection onto the quotient algebra. One can readily verify
that a necessary and sufficient condition that an operator family $A\in\cA$
is elliptic with a parameter is that the corresponding element $\si(A)$ of
the quotient algebra is invertible.  (As usual, when proving that the
family is invertible for large parameter values, one applies the Neumann
series.) Hence this condition can be used as an (equivalent) definition of
ellipticity:
\begin{definition}
The element $\si(A)$ is called the \textit{symbol} of a local operator
$A\in\cA$ with a parameter. The operator $A\in\cA$ is said to be
\textit{elliptic with a parameter} if its symbol $\si(A)$ is invertible.
\end{definition}
The equivalence of the two definitions is none other than the finiteness
theorem.
\begin{theorem}
If a family $A\in\cA$ is elliptic with a parameter, then it is Fredholm
with a parameter. The converse is also true.
\end{theorem}

Thus the main analytical task of elliptic theory of operators with a
parameter (just as in the parameter-free case) is the study of symbols
$\si(A)$ and the structure of the symbol algebra $\cA\slash\cJ$; in
particular, it is of interest to find specific subalgebras of $\cA$ for
which the symbols can be described constructively.

\section{Localization}

However, it turns out that under certain additional conditions there are
some nontrivial assertions even in elliptic theory for the \textit{entire}
algebra $\cA$. Namely, the invertibility of an element
$\si(A)\in\cA\slash\cJ$ is reduced to that of a set $\{p_x(A)\}$ of ``local
representatives'' labelled by points $x\in X$.

\subsection{Localizing classes}

Throughout the following, we assume that $X$ is a Hausdorff compactum. For
an arbitrary point $x\in X$, consider the set $\cF_x\subset C(X)$ of
functions $\ph(y)$ such that $0\le\ph(y)\le 1$ and $\ph(y)=1$ in some
(depending on $\ph$) neighborhood of $x$. This set is bounded (the norm of
each element is equal to~$1$) and multiplicative. It follows from the
Urysohn lemma that the ordering
\begin{equation*}
  \ph\prec\psi \Longleftrightarrow \ph\psi=\psi
\end{equation*}
makes it a directed set, which will be called the \textit{localizing class}
at $x$ (cf.~\cite[\S\,5.1]{GoKru2}). The supports of its elements
``shrink'' to $x$.

\subsection{Conditions on the action of $C(X)$}

We assume that the following two conditions are satisfied.
\begin{enumerate}
    \item[1$^\circ\!.$] The representation of $C(X)$ in $\cB H$ is faithful
                    (i.e., $C(X)$ is embedded in $\cB H$ as a subalgebra),
                    and $C(X)\cap \cK H=\{0\}$.
    \item[2$^\circ\!.$] For each $x\in X$, the localizing class
    $\cF_x$ strongly converges to zero in $\cB H$,
\begin{equation*}
    \slim_{\ph\in\cF_x}\ph =0.
\end{equation*}
\end{enumerate}
In other words, $H$ does not contain elements ``concentrated at $x$.'' (If
$H=L^2(X,d\mu)$, then condition 2$^\circ$ means that the measure $\mu$ is
purely nonatomic.)

\subsection{Localization in the algebra $\cA$}

We construct local representatives of elements $A\in\cA$ on the basis of
the localization principle in $C^*$-algebras, which we use in the form
given in~\cite[Proposition~3.1]{PlSe6}.

The algebra $C(X)$ is naturally embedded in $\cA$ as a commutative
subalgebra. (Its elements are constant functions of the parameter $y$). Let
$\cI_x\subset C(X)$ be the maximal ideal of functions vanishing at $x\in
X$, and let $\cJ_x$ be the ideal generated by $\cI_x$ in $\cA$. The
quotient algebra $\cA_x=\cA\slash\cJ_x$ is called the \textit{local
algebra} (at $x$), and the coset $p_x(A)\in\cA_x$ of an element $A\in\cA$
is called the \textit{local representative} of $A$ at $x$. Here
\begin{equation*}
    p_x:\cA\lra\cA_x
\end{equation*}
is the natural projection.

\begin{theorem}\label{loc-in-cA}
\rom{(a)} The $C^*$-algebra homomorphism
\begin{equation*}
    \left(\prod_{x\in X}p_x\right):\cA\lra\prod_{x\in X}\cA_x
\end{equation*}
induces a well-defined homomorphism
\begin{equation}\label{homo-mono}
\begin{aligned}
    \cA\slash\cJ&\lra\prod_{x\in X}\cA_x,\\
    \si(A)&\longmapsto\{p_x(A)\}_{x\in X}
\end{aligned}
\end{equation}
of the symbol algebra.

\rom{(b)} The homomorphism~\eqref{homo-mono} is a monomorphism\rom; that
is, a family $A\in\cA$ belongs to $\cJ$ \rom(i.e., is compact and decays at
infinity\rom) if and only if all local representatives $p_x(A)$, $x\in X$,
are zero.

\rom{(c)} A family $A\in\cA$ is elliptic  \rom(and hence Fredholm\rom) with
a parameter if and only if all local representatives $p_x(A)$, $x\in X$,
are invertible.
\end{theorem}

\begin{proof}
1. First, let us prove that the algebra $\cA$, the commutative subalgebra
$\cC=C(X)$, and the ideal $\cJ$ satisfy the assumptions of Proposition~3.1
in \cite{PlSe6}. Specifically, we should verify that
\begin{enumerate}
    \item[(i)] the restriction to $\cC$ of an arbitrary irreducible representation
      of $\cA$ is nonzero;
    \item[(ii)] for two arbitrary distinct points $x_1,x_2\in X$, there exist elements
    $\ph_1,\ph_2\in C(X)$ such that $\ph_1(x_1)\ne0$, $\ph_2(x_2)\ne0$, and
    $\ph_1\cA\ph_2\subset\cJ$;
    \item[(iii)] for an arbitrary irreducible representation $\pi$ of $\cJ$
    and any $x\in X$, there exist elements $A\in\cJ$ and $\psi\in\cI_x$
    such that $\pi(A\psi)\ne0$.
\end{enumerate}

Condition (i) is trivial, since the subalgebra $\cC$ contains the unit.
Condition (ii) holds, since for $\ph_{1,2}$ we can take functions supported
in disjoint neighborhoods of $x_{1,2}$; then the desired inclusion follows
from~\eqref{commut}.

To prove (iii), note that each irreducible representation $\pi$ of $\cJ$ is
equivalent to the representation $\pi(A)=A(y)\in\cB H$ for some
$y=y(\pi)\in Y$~\cite[Corollary~10.4.4]{Dix1}.

\medskip

Now take an element $A\in \cJ$ such that $A(y)\ne0$ and a vector $u\in H$
such that $A(y)u\ne0$. Since $\lim_{\ph\in\cF_x}(1-\ph)u=u$
(assumption~2$^\circ$), it follows that $A(y)(1-\ph)u\ne0$ for some
$\ph\in\cF_x$, and it suffices to set $\psi=1-\ph\in\cI_x$.

\medskip

2. Thus the assumptions of Proposition 3.1 in \cite{PlSe6} are satisfied.
The argument in the proof of this proposition shows that
\begin{equation}\label{loc-in-cA-eq}
    \wh{\cA\slash\cJ}=\bigcup_{x\in X}\wh\cA_x,
\end{equation}
where $\wh\cB$ is the \textit{spectrum} of the $C^*$-algebra $\cB$, i.e.,
the set of equivalence classes of its irreducible representations.

\medskip

Relation~\eqref{loc-in-cA-eq} is understood as follows. Irreducible
representations of the quotient algebras $\cA\slash\cJ$ and
$\cA_x=\cA\slash\cJ_x$ can be treated as irreducible representations of the
algebra $\cA$ itself. Then the right- and left-hand sides of the relation
are sets of equivalence classes of irreducible representations of $\cA$,
and the assertion is that these two sets coincide.

\medskip

It follows from~\eqref{loc-in-cA-eq} that
\begin{equation}\label{vazhno}
    \cJ=\bigcap_{x\in X}\cJ_x.
\end{equation}
Indeed, by \cite[Corollary, p.~34]{Arv1} for each element $A\notin\cJ$
there exists an irreducible representation $\pi\in\wh{\cA\slash\cJ}$ such
that $\pi(A)\ne 0$. It follows from~\eqref{loc-in-cA-eq} that
$\pi\in\wh{\cA\slash\cJ_x}$ for some $x$, so that $\ker\pi\supseteq\cJ_x$
and $A\notin\cJ_x$. Conversely, if $A\notin\cJ_x$ for some $x$, then,
reversing the argument, we obtain $A\notin\cJ$.

\medskip

In turn,~\eqref{vazhno} readily implies  (a) and (b).

3. Finally, to obtain (c),  on both sides of~\eqref{loc-in-cA-eq} we use
the fact that an element of a $C^*$-algebra is invertible if and only if so
are its images under all irreducible representations of the algebra (e.g.,
see \cite[Proposition~5.3]{PlSe7}).
\end{proof}

\subsection{Properties of local representatives}

One can naturally ask how the local representatives of elements $A\in\cA$
can be described more precisely.

First, we give an explicit description of the ideal $\cJ_x$ (and hence of
the quotient algebra $\cA_x$). It turns out that two elements of $\cA$ give
rise to the same element of $\cA_x$ if and only if they are equivalent with
respect to the localizing class $\cF_x$ in the sense of
formula~\eqref{e-opi1} below (cf.\ the definition
in~\cite[\S\,5.1]{GoKru2}).

\begin{proposition}\label{p-opis}
The ideal $\cJ_x$ consists of elements $A\in\cA$ such that
\begin{equation}\label{e-opi1}
    \lim_{\ph\in\cF_x}\bnorm{A\ph}=0.
\end{equation}
\end{proposition}
\begin{proof}
If an element $A\in\cA$ satisfies~\eqref{e-opi1}, then
\begin{equation*}
    A=\lim_{\ph\in\cF_x} A(1-\ph),
\end{equation*}
and since $1-\ph\in\cI_x$, we see that $A\in\cJ_x$. To prove the converse,
note that $\cJ_x$ is by definition the closure of the set of elements
representable as finite sums
\begin{equation}\label{finite-sum}
    C=\sum_nA_n\ph_nB_n,\quad A_n,B_n\in\cA,\quad \ph_n\in\cI_x.
\end{equation}
Commuting $B_n$ with $\ph_n$ and observing that each $\ph_n$ can be
approximated in norm by elements of the form $\ph_n(1-\psi)$,
$\psi\in\cF_x$, we see that $C$ can be approximated in norm by elements of
the form
\begin{equation}\label{odin-chlen}
    A=\wt C(1-\psi)+K,\quad \wt C =\sum A_nB_n\ph_n \in\cA,\quad K\in\cJ.
\end{equation}
Thus the elements~\eqref{odin-chlen} are dense in $\cJ_x$, and it suffices
to prove~\eqref{e-opi1} for these elements. In this case,
\begin{equation*}
    \lim_{\ph\in\cF_x} A\ph = \lim_{\ph\in\cF_x}\wt C(1-\psi)\ph
    + \lim_{\ph\in\cF_x}K\ph.
\end{equation*}
The first term is eventually zero, and it remains to prove that the limit
of the second term is zero. Since $K=K(y)\to 0$ as $y\to\infty$, it follows
that for each $\e>0$ there exist finitely many compact operators
$K_0=0,K_1,\dotsc,K_N$, $N=N(\e)$, such that for each $y\in Y$ the relation
\begin{equation}\label{priblizh}
    \norm{K(y)-K_j}\le\e
\end{equation}
is valid for  some $j=j(y)\in\{0,1,\dotsc,N\}$. It follows from condition
2$^\circ$ and the compactness of $K_j$ that
\begin{equation*}
    \lim_{\ph\in\cF_x}K_j\ph=0.
\end{equation*}
Combining this with~\eqref{priblizh}, we see that
\begin{equation*}
    \bnorm{\lim_{\ph\in\cF_x} K\ph}=\sup_y\bl\|\lim_{\ph\in\cF_x} K(y)\ph\br\|\le\e,
\end{equation*}
and we arrive at the desired assertion, since $\e$ is arbitrary.
\end{proof}

\begin{remark}
In~\eqref{e-opi1}, $\bnorm{A\ph}$ can be replaced with $\bnorm{\ph A}$ or
even with the norm of the element $\si(A\ph)=\si(\ph A)$ in $\cA\slash\cJ$.
\end{remark}

\medskip

Note that Theorem~\ref{loc-in-cA} does not describe the range of the
homomorphism~\eqref{homo-mono}. General facts about homomorphisms of
$C^*$-algebras only imply that this is a $C^*$-subalgebra of $\prod_x\cA_x$
isomorphic to $\cA\slash\cJ$. The following questions are natural to ask.
Let $\{a_x\}_{x\in X}$ be a family of elements of the local algebras
$\cA_x$. \textit{Under what conditions can this family be obtained by the
localization of some element of $\cA$}? \textit{Knowing the family, how can
one reconstruct this element modulo the ideal $\cJ$}?

The family $\{a_x\}$ can be viewed as a section of the projection
$\bigsqcup_x\cA_x\lra X$, where $\bigsqcup_x\cA_x$ is the disjoint union of
the sets $\cA_x$. The topology on $\bigsqcup_x\cA_x$ with respect to which
the families obtained by localization are exactly the \textit{continuous}
sections of the projection was studied in the framework of general
localization principles in $C^*$-algebras in many papers (e.g.,
see~\cite{DaHo1,Vas3} etc.). In our case, one can readily describe such
families without using the above-mentioned general results.

We introduce the following notation. Let $Q\subset X$ be an arbitrary
subset, and let $B\in C(Y,\cB H)$. By $\bnorm{B}_Q$ we denote the norm of
the restriction of $B$ to the subset $H_Q\subset H$ of elements supported
in the closure of $Q$:
\begin{equation*}
    \bnorm{B}_Q=\sup_{y\in Y} \norm{B(y)\bigm|_{H_Q}:H_Q\lra H}.
\end{equation*}

For each element $a_x$, we take some representative $A_x\in a_x$.
\begin{definition}\label{de1}
The family $\{a_x\}$ is said to be \textit{continuous} if for each $\e>0$
every point $x\in X$ has a neighborhood $U(\e,x)$ such that the following
condition is satisfied:
\begin{equation}\label{e-usl}
    \bnorm{A_x-A_{x'}}_{U(\e,x)\cap U(\e,x')}\le\e\quad \text{for any $x,x'\in X$.}
\end{equation}
\end{definition}
It is easily seen that this condition is independent of the choice of
representatives $A_x\in a_x$ (although the neighborhoods $U(\e,x)$ may
depend on this choice).

\medskip

Indeed, it follows from Proposition~\ref{p-opis} that if $\wt A_x$ are some
other representatives of $a_x$, then for each $x$ then norm $\bnorm{A_x-\wt
A_x}_U$ can be made as small as desired if we take a sufficiently small
neighborhood $U$ of the point $x$ (where the corresponding element
$\ph\in\cF_x$ is equal to unity).

\begin{proposition}\label{p2}
The family $\{a_x\}$ is the localization of some operator $A\in\cA$ if and
only if it is continuous.
\end{proposition}
\begin{proof}
Let $a_x=\si_x(A)$, $x\in X$. Then one can take $A_x=A$ and $U(\e,x)=X$ for
all $x\in X$ and $\e>0$. Then, obviously, $\bnorm{A_x-A_{x'}}=0$, so that
condition~\eqref{e-usl} holds.

Conversely, if $\{a_x\}$ is a continuous family of local representatives,
then an element $A\in\cA$ such that
\begin{equation}\label{e1}
a_x=p_x(A)\quad\text{for all $x\in X$},
\end{equation}
can be constructed by the following method. For each element $a_x$, take
some representative $A_x\in a_x$. Next, for each $\e>0$ take a finite
continuous nonnegative partition of unity $1=\sum_{x\in X}\ph_{\e x}$ on
$X$ subordinate to the cover $\{U(\e,x)\}$ in Definition~\ref{de1} and set
\begin{equation*}
    A^{(\e)}=\sum_x \ph_{\e x}A_x.
\end{equation*}
We claim that there exists a limit (independent of the ambiguity in the
construction)
\begin{equation*}
    B=\lim_{\e\to0}\si\bl(A^{(\e)}\br)\in\cA\slash\cJ
\end{equation*}
and that~\eqref{e1} holds for arbitrary $A\in\cA$ such that $\si(A)=B$.
Indeed,
\begin{equation}\label{zvzda}
    A^{(\e)}-A^{(\dt)}=\sum_{x,y}\ph_{\e x}\ph_{\dt y}(A_x-A_y).
\end{equation}
Let us use the following lemma.
\begin{lemma}\label{lemma8}
Let $f_j\in C(X)$, $j=1,\ldots,N$, be nonnegative function, and let
$A_j\in\cA$ be some elements. Then
\begin{equation}\label{e-chudo}
    \Bl\|\si\Bl(\sum_{j=1}^Nf_jA_j\Br)\Br\|_{\cA\slash\cJ}
    \le\Bl[\max_{x\in X}\sum_{j=1}^N f_j(x)\Br]\max_{j=1,\dots,N}
    \bnorm{A_j}_{\supp f_j}.
\end{equation}
\end{lemma}
The proof of this lemma is similar to that of Theorem~4.1 in~\cite{GoKru2}.
It follows from the lemma that
\begin{equation*}
    \norm{\si\bl(A^{(\e)}\br)-\si\bl(A^{(\dt)}\br)}_{\cA\slash\cJ}\le\max\{2\e,2\dt\},
\end{equation*}
so that $\si\bl(A^{(\e)}\br)$ is a Cauchy sequence in $\cA\slash\cJ$ and
hence has a limit, $B$. Using relations like~\eqref{zvzda}, one can readily
show that this limit is independent of the choice of partitions of unity.

Let $\si(A)=B$; let us show that~\eqref{e1} holds. We take a point $x_0\in
X$ and choose partitions of unity such that $\ph_{\e x_0}(x)=1$ in some
neighborhood (depending on $\e$) of $x_0$. Then, obviously,
$p_{x_0}A^{(\e)}=a_{x_0}$. Moreover,
\begin{equation*}
    A=\lim_{\e\to0} \bl(A^{(\e)}+K^{(\e)}\br),
\end{equation*}
where the $K^{(\e)}$ are appropriate elements of the ideal $\cJ$. Then
\begin{equation*}
    p_{x_0}(A)=\lim_{\e\to0} p_{x_0}\bl(A^{(\e)}+K^{(\e)}\br)
            =\lim_{\e\to0} p_{x_0}\bl(A^{(\e)}\br)=a_{x_0}.
\end{equation*}
The proof of the proposition is complete.
\end{proof}

Summarizing, we obtain the following theorem.
\begin{theorem}\label{t0}
The set $\Sigma\subset\bigsqcup_{x\in X}\cA_x$ of continuous families is a
$C^*$-algebra with respect to the norm
\begin{equation*}
    \norm{\{a_x\}}=\sup_{x}\norm{a_x}_{\cA_x}.
\end{equation*}
The mapping $\si(A)\longmapsto\{p_x(A)\}$ is a well-defined isometric
isomorphism of the algebra $\cA\slash\cJ$ onto the algebra $\Sigma$.
\end{theorem}

\section{Cones and infinitesimal operators}

An efficient application of Theorem~\ref{t0} to general local operators is
hindered by the fact that the local algebras $\cA_x$ are not naturally
realized as operator algebras. In particular, one cannot choose a
``canonical'' representative $A_x$ of $a_x$ so as to verify the
invertibility of $a_x$ conveniently. However, if $X$ bears some additional
structures, then in $\cA$ there exist subalgebras for which the
corresponding local algebras are naturally realized as operator algebras.
In this section, we describe a relevant construction for the case in which
the additional structure is the ``tangent cone'' to $X$. For the case of
operators without parameter, a close construction can be found
in~\cite[Chap. 4]{AnLe1}.

\subsection{Operators on the cone}

By a \textit{cone} we mean a noncompact locally compact Hausdorff space $K$
with a distinguished point $0\in K$ (the vertex) and with continuous
one-parameter multiplicative dilation group
\begin{equation*}
    g_\la:K\lra K,\quad\la\in\RR_+,
\end{equation*}
which leaves the vertex fixed ($g_\la(0)=0$) and has the following
property: for each compact set $Q\subset K$ and each neighborhood of $V$,
\begin{equation*}
 g_\la(Q)\subset V\quad \text{for all sufficiently small $\la\in\RR_+$.}
\end{equation*}
it follows that the family of functions
\begin{equation}\label{phila}
    \ph_\la(x)=\ph(g_\la x),\quad x\in K,\quad\la\to\infty,
\end{equation}
is cofinal in $\cF_0$ for each $\ph\in\cF_0$, where $\cF_0$ is the set of
compactly supported continuous functions $\ph(x)$ on $K$ such that
$0\le\ph(x)\le 1$ and $\ph(x)=1$ in some neighborhood (depending on $\ph$)
of the vertex.

Next, let a Hilbert space $\wt H$ be a $*$-module over the $C^*$-algebra
$C_0(K)$ of continuous functions on $K$ decaying at infinity, and suppose
that the strong closure of the image of $C_0(K)$ in $\cB\wt H$ contains the
identity operator and assumption $2^\circ$ holds. Finally, let
\begin{equation*}
 U_\la:\wt H\lra \wt H, \quad\la\in\RR_+,\qquad U_\la U_\mu=U_{\la\mu},
\end{equation*}
be a strongly continuous one-parameter group of unitary operators in $\wt
H$ related to the action of $C_0(K)$ by the formula
\begin{equation*}
    U_\la \ph(x) U_\la^{-1}=\ph(g_\la x),\quad\la\in\RR_+,\quad x\in K.
\end{equation*}
\begin{proposition}
The group $U_\la$ weakly converges to zero in $\wt H$ as $\la\to0$ as well
as as $\la\to\infty$.
\end{proposition}
\begin{proof}
It suffices to prove the weak convergence to zero as $\la\to\infty$; the
second assertion then follows, since the group $U_\la$ is unitary. Next, by
the Banach--Steinhaus theorem it suffices to prove the weak convergence on
the dense subset of elements of the form $\ph u$, where $u\in \wt H$ and
$\ph\in\cF_0$. We have
\begin{equation*}
    (v, U_\la\ph u)=(v, U_\la\ph U_\la^{-1}U_\la u)
                   =(v,\ph_\la U_\la u)=(\ph_\la v,U_\la u)\to 0
                   \quad \text{as $\la\to\infty$,}
\end{equation*}
since $\ph_\la v\ovs{s}\lra 0$ by assumption $2^\circ$ ($\ph_\la$ is
cofinal in $\cF_0$) and the family $U_\la u$ is uniformly bounded.
\end{proof}

Suppose that $\RR_+$ also acts by homeomorphisms on the parameter space
$Y$. This action will be denoted by $(\la,y)\longmapsto \la y$.

By $\cL\subset C(Y,\cB\wt H)$ we denote the closed $C^*$-subalgebra of
families $B$ such that
\begin{equation*}
 [B,\ph]\in\cJ=C_0(Y, \cK\wt H)\quad\text{for each $\ph\in C_0(K)$}.
\end{equation*}
Next, $\cJ_0\subset\cL$ is the closed $*$-ideal in $\cL$ consisting of the
families $B\in\cL$ such that $B\ph\in\cJ$ for each $\ph\in
C_0(K)$.\footnote{It obviously suffices to impose this condition for
$\ph\in\cF_0$.} Elements of this ideal will be called \textit{almost
compact operators with a parameter}. An element $B\in\cL$ is said to be
\textit{homogeneous} if
\begin{equation}\label{eq-odn}
 U_\la^{-1}B(\la y)U_\la=B(y)\quad\text{for all $\la\in\RR_+$}.
\end{equation}
The set of all homogeneous elements will be denoted by $\cL_\infty$.
Obviously, this is a $C^*$-subalgebra in $\cL$.

\begin{proposition}
\begin{equation*}
    \cL_\infty\cap\cJ_0=\{0\}.
\end{equation*}
\end{proposition}
\begin{proof}
Let $B\in \cL_\infty\cap\cJ_0$. Consider the family $B\ph$, where
$\ph\in\cF_0$. Obviously, if $\psi\in\cF_0$, then for $\la\gg1$ we have
$\psi\ph_\la=\ph_\la$ and
\begin{equation*}
    \bnorm{B\ph}=\bnorm{U_\la B\ph U_\la^{-1}}
                =\sup_y\norm{B(\la y)\psi\ph_\la}.
\end{equation*}
The family $B(y)\psi$ is contained in the ideal $\cJ$ and hence can be
approximated with arbitrary accuracy by a step function of $y$ assuming
finitely many compact values. Since $\ph_\la\ovs{s}\to0$, we conclude that
the right-hand side of the last relation tends to zero as $\la\to\infty$,
and hence $B\ph=0$ for each $\ph\in\cF_0$. Taking a sequence of functions
$\ph$ strongly convergent to the identity operator, we obtain $B=0$, as
desired.
\end{proof}

\subsection{Infinitesimal operators}

Now let $x\in X$, and let a homeomorphism
\begin{equation*}
    f: \ov V\lra f(\ov V)\subset K
\end{equation*}
of the closure of some neighborhood $V$ of $x$ on the closure of a
neighborhood of the vertex of $K$ be given. Next, suppose that an
isomorphism
\begin{equation*}
    \ga:H_{f(\ov V)}\lra \wt H_{\ov V}
\end{equation*}
of $*$-modules over the ring homomorphism $f^*:C(f(\ov V))\lra C(\ov V)$ is
given. (In what follows, we identify $V$ and $f(V)$ as well as the
subspaces related by $\ga$.) For brevity, we refer to $K$ as the
\textit{tangent cone} to $X$ at $x$. (Needless to say, we do not claim
uniqueness.)

\begin{definition}\label{inf-op}
An operator $A\in\cA$ with a parameter is said to be
\textit{infinitesimally localizable} at the point $x$ along the tangent
cone $K$ if the class $p_{x}(A)\in\cA_x$ contains a representative of the
form $\psi B\ph$, where $B\in\cL_\infty$ and $\psi,\ph\in\cF_x$ are
function supported in $V$. The operator $B$ is called the
\textit{infinitesimal operator} for $A$ at $x$ and is denoted by
\begin{equation*}
  B=i_x(A).
\end{equation*}
\end{definition}

\begin{theorem}\label{th13}
Definition~\rom{\ref{inf-op}} is consistent. The mapping $A\longmapsto
i_x(A)$ is a $C^*$-algebra homomorphism and factors through the quotient
algebra $\cA_x$.
\end{theorem}
\begin{proof}
1. First of all, let us prove the consistency, i.e., the fact that the
element $i_x(A)\in\cL_\infty$ is uniquely determined. Indeed, let $B$ and
$\wt B$ be two such elements, so that
\begin{equation*}
    \psi B\ph, \wt\psi\wt B\wt\ph\in p_x(A).
\end{equation*}
Set $\Delta=B-\wt B\in\cL_\infty$. Taking a function $\chi\in\cF_x$ such
that all four elements of $\cF_x$ in the last formula are equal to unity on
the support of $\chi$, we obtain
\begin{equation*}
    \chi\Delta\chi=\chi(B-\wt B)\chi=\chi(\psi B\ph-\wt\psi\wt
    B\wt\ph)\chi\in\cJ_x.
\end{equation*}
By applying Proposition~\ref{p-opis}, we obtain
\begin{equation*}
    \lim_{\rho\in\cF_x}\bnorm{\chi\Delta\chi\rho}=0.
\end{equation*}
Take $\rho=\phi_\la$, where $\phi\in\cF_0$. For sufficiently large $\la$,
we have $\chi\phi_\la=\phi_\la$ and
\begin{equation*}
    \bnorm{\phi_\la\Delta\phi_\la}\le\bnorm{\chi\Delta\chi\phi_\la}\lra0.
\end{equation*}
On the other hand, by passing to the cone and by using the homogeneity of
the element $\Delta$, we obtain
\begin{equation}\label{zzvezdochka}
    \begin{aligned}
\bnorm{\phi_\la\Delta\phi_\la}&=\bnorm{U_\la^{-1}\phi_\la\Delta\phi_\la
U_\la}
     =\bnorm{\phi U_\la^{-1}\Delta U_\la\phi}
     \\&=\sup_{y}\norm{\phi \Delta(y/\la)\phi}
     =\sup_{y}\norm{\phi \Delta(y)\phi}=\bnorm{\phi\Delta\phi}.
\end{aligned}
\end{equation}
Hence $\phi\Delta\phi=0$ for each $\phi\in\cF_0$, and consequently,
\begin{equation*}
    \Delta=\slim_{\la\to0}\phi_\la\Delta\phi_\la=0.
\end{equation*}
Now it is obvious that if $A\in\cJ_x$, then $i_x(A)=0$. Indeed, $0\in
p_x(A)$ in this case.

2. Let us verify that the mapping $i_x$ is multiplicative. (Additivity is
obvious.) Let $A$ and $C$ be infinitesimally localizable operators. Then
\begin{multline*}
    AC-\psi i_x(A)i_x(C)\ph=AC-\psi i_x(A)\chi^2i_x(C)\ph+K\\
    =(A-\psi i_x(A)\chi)C+\psi i_x(A)\chi(C-\chi i_x(C)\ph)+K
\end{multline*}
for some $\chi\in\cF_x$ and $K\in\cJ$, and all terms on the right-hand side
belong to $\cJ_x$. Thus the product $AC$ is infinitesimally localizable,
and
\begin{equation*}
    i_x(AC)=i_x(A)i_x(C).
\end{equation*}
3. Let us prove that the set of operators infinitesimally localizable at
$x$ is closed. Using the same trick as in~\eqref{zzvezdochka}, we can show
that for a homogeneous element $i_x(A)\in\cL_\infty$ one always has
\begin{equation*}
    \bnorm{i_x(A)}=\bnorm{\ph i_x(A)\ph},\quad \ph\in\cF_0.
\end{equation*}
On the other hand,
\begin{equation*}
    \bnorm{\ph(A-i_x(A))\ph}\ovs{\ph\in\cF_x}\lra0.
\end{equation*}
It follows that
\begin{equation}\label{estimate1}
\bnorm{i_x(A)}\le\bnorm{A}.
\end{equation}
Since the set of homogeneous elements is a $C^*$-algebra and, in
particular, is closed, we see that so is the set of infinitesimally
localizable operators. The proof is complete.
\end{proof}

Let $Z$ be a set of pairs of the form $(x,K)$, where $x\in X$ and $K$ is a
tangent cone at $x$. We define a projection $\pi:Z\lra X$ by setting
$\pi(x,K)=x$. (We do not exclude the case in which several tangent cones
are given at some point $x$.) We assume that $\pi(Z)=X$. We denote by
$\cA_Z$ the subset of $\cA$ formed by the elements infinitesimally
localizable at $x$ along $K$ for each pair $(x,K)\in Z$. The preceding
argument readily implies the following theorem.
\begin{theorem}
The set $\cA_Z$ is a $C^*$-algebra. The mapping
\begin{equation*}
    A\longmapsto\{i_x(A)\}_{(x,K)\in Z}
\end{equation*}
factors through $\cA_Z\slash\cJ$ and is a monomorphism.
\end{theorem}

\subsection{The decomposition problem for infinitesimal operators}

If a local operator $A\in\cA$ with a parameter has an infinitesimal
operator $i_x(A)$ (along a given cone) at a point $x$ and if this
infinitesimal operator is invariant with respect to some transformation
group, then, passing to the Fourier transform associated with this group,
one can decompose the operator $i_x(A)$ into a direct sum (or direct
integral) over irreducible representations of the group. This was
indicated, say, in~\cite[Chap.~4]{AnLe1}. A typical example is given by
classical zero-order $\Psi$DO $A$ on a smooth manifold $X$. Here the
tangent cone is just the tangent space $T_xX$ at $x$ (a neighborhood of
zero in the tangent space is mapped onto a neighborhood of the point in the
manifold, say, with the help of local coordinates with origin at $x$), the
infinitesimal operator $i_x(A)$ is a translation invariant operator in
$L^2(T_xX)$ and hence has the form
\begin{equation*}
    i_x(A)=f\BL(x,-i\pd{}{t}\BR)
\end{equation*}
(where $x$ is a parameter and the variables $t$ are coordinates on $T_xX$),
and the Fourier transform associated with the translation group is none
other than the ordinary Fourier transform and takes the infinitesimal
operator to the operator of multiplication by a function $f(x,\xi)$,
$\xi\in T^*_xX$. This function is called the (principal) symbol of $A$, is
denoted by $\si(A)(x,\xi)$, is continuous and zero-order homogeneous with
respect to $\xi$ for $\xi\ne0$, and has a jump discontinuity at the zero
section of the cotangent bundle $T^*X$.

What are the conditions under which the infinitesimal operator proves to be
invariant under some transformation group, so that the Fourier transform
gives a ``symbol''? We give an answer in the next subsection.

\subsection{Symbols in the small}

Let $K$ be a cone (whose points will be denoted by $z$ and the action of
elements $\la\in\RR_+$ by $\la z$), and let $\wt H$ be a Hilbert space that
is a $*$-module over $C_0(K)$ and is equipped with a unitary group
associated with dilations. (Here the group will be denoted by $\ka_\la$.)
Next, suppose that an open subset $U\subset X$ is homeomorphically mapped
onto a neighborhood of the point $(0,0)$ in the Cartesian product
$\RR^k\times K$. (For elements of $U$, we use the coordinates $(t,z)$,
$t\in\RR^k$, $z\in K$.) Finally, let the subspace $H_U$ be identified via
an isomorphism respecting the module structure with the corresponding
subspace of the Hilbert space $\gH=L^2(\RR^k,\wt H)$ (on which the algebra
$C_0(\RR^k\times K)$ acts naturally, i.e., pointwise with respect to the
argument $t\in\RR^k$). Then for each point of the form $(x,0)\in U$ the
product $\RR^k\times K$ can be viewed as the tangent cone with vertex
$(x,0)$, the action of $\RR_+$ given by the formula
\begin{equation*}
    g_{x,\la}(t,z)=\bl(x+\la(t-x),\la z\br), \quad\la\in\RR_+,
\end{equation*}
and the one-parameter unitary group in $\gH$ given by the
formula\footnote{Here the factor $\la^{k/2}$ ensures that the group is
unitary.}
\begin{equation*}
    U_{x,\la}f(t)=\la^{k/2}(\ka_\la f)\bl(x+\la(t-x)\br).
\end{equation*}
A straightforward computation shows that the group interacts with the
action of $C_0(\RR^k\times K)$ by the desired formula
\begin{equation*}
    U_{x,\la}\ph(t,z)U_{x,\la}^{-1}=\ph(g_{x,\la}(t,z)).
\end{equation*}

We are interested in infinitesimal operators of a given operator $A$ with a
parameter at the points $(x,0)\in U$. For brevity, we denote them by
$a_x=i_x(A)$ instead of $i_{(x,0)}(A)$. They all act in the space $\gH$,
but the cone is different for each $x$. (The cone vertex and the dilation
groups depend on $x$). These infinitesimal operators are determined by the
homogeneity condition~\eqref{eq-odn} and the convergence
\begin{equation*}
    \lim_{\phi\in\cF_x}\bnorm{\psi(A-a_x)\phi}=0,
\end{equation*}
where $\psi\in C_0(U)$ is an arbitrary cutoff function. We consider only
$x$ close to zero; hence it can be assumed that $A$ has already been
multiplied on the left and on the right by cutoff functions and is well
defined in $\gH$, so that the cutoff factor $\psi$ in this condition can be
omitted. Next, we represent the function $\phi\in\cF_x$ in the convergence
condition in the form
\begin{equation*}
    \phi=T_{-x}\ph T_x,
\end{equation*}
where $\ph\in\cF_0$ and $T_x$ is the translation operator
\begin{equation*}
    T_xf(t,z)=f(t+x,z).
\end{equation*}
Then the condition acquires the form
\begin{equation}\label{uniform}
    \lim_{\ph\in\cF_0}\bnorm{(A-a_x)T_{-x}\ph T_x}=0,
\end{equation}
which will be used in what follows.\footnote{The factor $T_{-x}\ph T_x$ can
also be placed on the left; the results will be the same.}

\begin{theorem}\label{999}
\rom{(i)} If $A\in\cA$ is an operator with a parameter such that the
infinitesimal operators $a_x=i_x(A)$ exist at all points $(x,0)\in U$ and
the convergence in~\eqref{uniform} is locally uniform with respect to $x$,
then the infinitesimal operators $a_x$ commute with translations with
respect to $x$, i.e., satisfy
\begin{equation*}
    [a_x,T_\tau]=0\quad\text{for any $\tau$ and $x$,}
\end{equation*}
and depend on $x$ continuously in the norm $\bnorm{\cdot}$. The set of
operators satisfying these conditions is a $C^*$-subalgebra in $\cA$ and
contains the ideal~$\cJ$.

\rom{(ii)} Conversely, suppose that the infinitesimal operators $a_x$ for
an operator $A\in\cA$ exist, commute with translations, and continuously
depend on the parameter $x$. Then the convergence in~\eqref{uniform} is
locally uniform with respect to $x$.
\end{theorem}
\begin{proof}
First, let us prove (i).

\textbf{1.} Let us verify that the infinitesimal operators are translation
invariant. It suffices to do this for $x=0$. In~\eqref{uniform}, we replace
the arbitrary function $\ph$ by $\ph_\la=U_\la\ph U_\la^{-1}$
(see~\eqref{phila}), where $U_\la\equiv U_{0,\la}$ and the function
$\ph\in\cF_0$ is arbitrary but fixed. Then we speak of convergence locally
uniform with respect to $x$ as $\la\to\infty$. First, we write out some
relations between our groups, which can be verified by a straightforward
computation:
\begin{equation}\label{sootn}
    U_{x,\la}=T_{-x}U_{0,\la}T_x,\quad U_\la
    T_x=T_{x/\la}U_\la=T_xU_{x,\la}.
\end{equation}
Throughout the following, we assume that $x$ varies in a neighborhood of
zero; moreover, $\la>1$, so that $x/\la$ lies in the same neighborhood.
Since condition~\eqref{uniform} is uniform with respect to $x$, we can
replace $x$ by $x/\la$ in it. Separately substituting also $x=0$ into this
condition, we obtain (the arrow means convergence in the norm
$\bnorm{\bcdot}$ as $\la\to\infty$)
\begin{equation}\label{1}
    (A-a_0)\ph_\la\to0,\quad (A-a_{x/\la})T_{-x/\la}\ph_\la
    T_{x/\la}\to0.
\end{equation}
We multiply the equations in~\eqref{1} on the left and on the right by the
unitary operators $T_{\pm x/\la}$ and obtain
\begin{equation}\label{2}
    T_{x/\la}(A-a_0)T_{-x/\la}T_{x/\la}\ph_\la T_{-x/\la}\to0,
    \quad T_{x/\la}(A-a_{x/\la})T_{-x/\la}\ph_\la
    \to0.
\end{equation}
Take a function $\chi\in\cF_0$ such that
\begin{equation*}
    \chi\ph=\chi T_{x}\ph  T_{-x}=\chi
\end{equation*}
for all $x$ in our neighborhood. This can be done if the space is
sufficiently small. The function $\chi$ can be assumed to be an arbitrary
given element of $\cF_0$. (Given $\chi$, there always exists an appropriate
$\ph$.) Multiplying both equations in~\eqref{2} on the right by $\chi_\la$
and using~\eqref{sootn}, we obtain
\begin{equation*}
    T_{x/\la}(A-a_0)T_{-x/\la}\chi_\la\to0,
    \quad T_{x/\la}(A-a_{x/\la})T_{-x/\la}\chi_\la
    \to0,
\end{equation*}
and then, subtracting one from the other,
\begin{equation}\label{3}
    T_{x/\la}(a_{x/\la}-a_0)T_{-x/\la}\chi_\la
    \to0.
\end{equation}
In~\eqref{3}, we use the homogeneity property
\begin{equation*}
    a_0(y)=U_\la a_0(y/\la)U_\la^{-1},\quad
    a_{x/\la}(y)=U_{x/\la,\la} a_{x/\la}(y/\la)U_{x/\la,\la}^{-1}
\end{equation*}
of infinitesimal operators. By substituting this into the left-hand side
of~\eqref{3} and by using relations~\eqref{sootn}, we obtain
\begin{equation*}
 \begin{split}
    T_{x/\la}(a_{x/\la}(y)-a_0(y))&T_{-x/\la}\chi_\la\\
    &=\bl[U_\la T_{x/\la}a_{x/\la}(y/\la)T_{-x/\la}U_\la^{-1}
    -T_{x/\la}U_\la a_{0}(y/\la)U_\la^{-1}T_{-x/\la}\br]\chi_\la
    \\
    & =U_\la\bl[T_{x/\la}a_{x/\la}(y/\la)T_{-x/\la}
    -T_{x}a_{0}(y/\la)T_{-x}\br]\chi U_\la^{-1}\to0.
\end{split}
\end{equation*}
The left- and rightmost unitary factors $U_\la^{\pm1}$ can be omitted.
Next, since the norm $\bnorm{\bcdot}$ involves the supremum over $y$, we
can replace the argument $y/\la$ of infinitesimal families by $y$, thus
obtaining
\begin{equation}\label{4}
 [T_{x/\la}a_{x/\la}(y)T_{-x/\la}
    -T_{x}a_{0}(y)T_{-x}\br]\chi\to0.
\end{equation}
It follows from~\eqref{4} that the operator
$T_{x/\la}a_{x/\la}(y)T_{-x/\la}\chi$ converges in the norm to the operator
$T_{x}a_{0}(y)T_{-x}\chi$ as $\la\to\infty$. On the other hand, one can
readily show that under the conditions of the theorem it strongly converges
to the operator $T_0a_0(y)T_0\chi=a_0(y)\chi$.

Indeed, let us show that the operator function
\begin{equation*}
    \wt a_x=T_xa_xT_{-x}\chi
\end{equation*}
is strongly continuous. It follows from the definitions that
\begin{equation*}
    \lim_{\la\to\infty} F(x,\la)=\wt a_x
\end{equation*}
locally uniformly with respect to $x$ in the norm and hence in the strong
sense, where the operator function
\begin{equation*}
    F(x,\la)=U_\la^{-1}T_xAT_{-x}U_\la\chi
\end{equation*}
is strongly continuous in $x$. Thus the desired assertion follows from the
continuity of a uniform limit of continuous functions.

The limits should coincide, and hence
\begin{equation}\label{5}
    T_{x}a_{0}(y)T_{-x}\chi=a_0(y)\chi
\end{equation}
for each $\chi\in\cF_0$. Since vectors of the form $\chi u$ are dense in
$\gH$ (recall that the identity operator lies in the strong closure of the
algebra $C_0(\RR^k\times K)$), we obtain the desired assertion
\begin{equation}\label{6}
    T_{x}a_{0}T_{-x}=a_0.
\end{equation}

\textbf{2.} Let us prove that the operator $a_x$ continuously depends on
$x$ for $x=0$. For a given $\e>0$, for each point $(x,0)\in U$ we take a
neighborhood $U(\e,x)$ in accordance with Proposition~\ref{p2}. We consider
only points $x$ that are so close to zero that $(x,0)\in U(\e,0)$; then the
intersection $O=U(\e,0)\cap U(\e,x)$ is necessarily a nonempty open set (a
neighborhood of the point $(x,0)$). By Proposition~\ref{p2},
$\bnorm{\psi(a_x-a_0)}_O<\e$, where $\psi\in\cF_0$ is the function equal to
unity in $O$. Since the operator $a_x-a_0$ is translation invariant, we
obtain
\begin{equation*}
    \bnorm{\wt\psi(a_x-a_0)}_{\wt O}<\e,
\end{equation*}
where $\wt O=T_{-x}(O)$ is a neighborhood of $(0,0)$ and the translate
$\wt\psi$ of $\psi$ by $-x$ again lies in $\cF_0$. Multiplying this on the
right by a function $\ph\in\cF_0$ supported in $\wt O$, we obtain
\begin{equation*}
    \bnorm{\wt\psi(a_x-a_0)\ph}_{\wt O}<\e.
\end{equation*}
The operator $B=a_x-a_0$ has the homogeneity property~\eqref{eq-odn} with
respect to the group $U_\la$. (For $a_0$ this is true by definition, and
for $a_x$ this follows from the definition combined with the translation
invariance.) Using the homogeneity, we make the supports of $\wt\psi$ and
$\ph$ arbitrarily large without changing the operator $B$. It remains to
use the fact that the identity operator belongs to the strong closure of
$C_0(\RR^k\times K)$ to conclude that $\bnorm{a_x-a_0}<\e$ for $x$ close to
zero.

\medskip

\textbf{3.} Now let us show that the elements of $\cJ$ satisfy the
conditions of the theorem. For an element $A\in\cJ$,
condition~\eqref{uniform} can be rewritten as
\begin{equation}\label{uniform-comp}
    \lim_{\ph\in\cF_0}\bnorm{AT_{-x}\ph}=0.
\end{equation}
(We have omitted the unitary factor $T_x$ that does not affect the norm and
used the fact that $a_x=0$ for a compact operator.) Since $A$ is compact
and the operator $T_x$, which is the adjoint of $T_{-x}$, is a strongly
continuous function of $x$, it follows that the family $A_x=AT_{-x}$ is
continuous in the norm and condition, and so the convergence
in~\eqref{uniform-comp} is locally uniform in $x$.

\textbf{4.} Let us show that the set of elements satisfying the assumptions
of item~(i) is a $*$-algebra. Assertion (ii) of the theorem, proved below,
permits us to replace the assumptions of item (i) by those of item (ii),
which obviously are respected by algebraic operations and the passage to
the adjoint operator.

\medskip

It remains to show that this set is closed. Let a sequence $A_n\in\cA$
converge in $A$, and suppose that for each $n$ the convergence
\begin{equation*}
        \lim_{\ph\in\cF_0}\bnorm{(A_n-a_{nx})T_{-x}\ph T_x}=0
\end{equation*}
is locally uniform in $x$. As was shown earlier (Theorem~\ref{th13}), the
limit element $a$ is also infinitesimally localizable at each point
$(x,0)$, and moreover, $a_x=\lim_{n\to\infty} a_{nx}$. The convergence of
$a_{nx}$ to $a_x$ is uniform, since (see~\eqref{estimate1})
\begin{equation*}
 \bnorm{a_n-a_{nx}}\le\bnorm{A-A_n}.
\end{equation*}
Now it is obvious that the convergence in~\eqref{uniform} for the limit
operator $A$ is uniform.

\bigskip

Now let us prove (ii). We give only a brief scheme of the proof; the reader
can readily reconstruct the details. We use the method in
Proposition~\ref{p2} to reconstruct the operator $A$ from its local
representatives. Moreover, we choose the neighborhoods $U(\e,(\xi,z))$ of
points$(\xi,z)$ with $z\ne0$ (i.e., $z$ is not the cone vertex) so that
they do not contain any points of the form $(x,0)$. For the local
representatives we take the elements $a_x$ for the points $(x,0)$ and $A$
for all other points. Having constructed the approximations by this method,
we have
\begin{equation*}
    A=\lim_{\e\to0} \bl(A^{(\e)}+K^{(\e)}\br),
\end{equation*}
where the $K^{(\e)}$ belong to $\cJ$. Using the assumptions imposed on the
local representatives, we can readily verify that the operators $A^{(\e)}$
(and hence $A^{(\e)}+K^{(\e)}$) satisfy all assumptions of item (i) of the
theorem; namely, they are infinitesimally localizable at each point
$(x,0)$, and the corresponding convergence~\eqref{uniform} is locally
uniform. But then, as was shown above, the uniform convergence holds also
for the limit operator $A$.

The proof of the theorem is complete.
\end{proof}

Thus the theorem gives the desired group commuting with infinitesimal
operators. Now we can use the Fourier transform and decompose the
infinitesimal operators into families of simpler operators depending on the
parameter $\xi$. The following assertion establishes the corresponding
result.

\begin{proposition}\label{auxp1}
Let $H$ be a \rom(separable\rom) Hilbert space, and let
\begin{equation*}
 \wh B:L^2(\RR^k,H)\lra L^2(\RR^k,H)
\end{equation*}
be a bounded linear operator commuting with translations in $\RR^k$\rom:
\begin{equation*}
    \wh BT_x=T_x\wh B, \quad x\in\RR^k,\qquad\text{where}\quad [T_xu](t)=u(x+t).
\end{equation*}
Then $B$ can be represented in the form
\begin{equation*}
    \wh B=\ov\cF_{\xi\to t} B(\xi) \cF_{t\to\xi},
\end{equation*}
where $\cF$ is the Fourier transform and $B(\xi)$ is a strongly
measurable\footnote{In the separable case, strong and weak measurability
are equivalent by the Pettis theorem~\cite{Yos1}.} operator-valued function
ranging in $\cB H$ and satisfying the estimate
\begin{equation*}
    \esssup_\xi\norm{B(\xi)}=\norm{\wh B}.
\end{equation*}
If the operator $\wh B$ continuously depends on some parameters $\tau$,
then $B(\xi)$ continuously depends on the same parameters\rom:
\begin{equation*}
    \esssup_\xi\norm{B_\tau(\xi)-B_{\tau_0}(\xi)}\to 0\quad \text{as $\tau\to \tau_0$.}
\end{equation*}
\end{proposition}

The proof will be given in the Appendix.

\begin{remark}
If a function $u\in L^2(\RR^k,H)$ is such that its Fourier transform $\wt
u$ lies in $L^2(\RR^k,H)\cap L^1(\RR^k,H)$, then $B\wt u$ lies in the same
space, and hence the inverse Fourier transform is given by a Bochner
integral and we can write
\begin{equation*}
    [Bu](t)=\BL(\frac{1}{2\pi}\BR)^{n/2}\int e^{it\xi}B(\xi)\wt
    u(\xi)\,d\xi.
\end{equation*}
\end{remark}

By the proposition we have just proved, to the infinitesimal operator
$i_x(A)$ there corresponds an operator-valued function $\si(x,\xi,y)$ such
that
\begin{equation}\label{sootv}
   i_x(A)(y)=\ov\cF_{\xi\to t} \si(x,\xi,y) \cF_{t\to\xi},
\end{equation}
(here $(x,y)$ play the role of additional parameters $\tau$). Using the
definition of the group $U_\la$ and the homogeneity property $i_x(A)(\la
y)=U_\la i_x(A)(y)U_\la^{-1}$ and passing to the Fourier transform, we see
that this operator-valued function satisfies the \textit{twisted
homogeneity} condition
\begin{equation}\label{skruch}
    \si(x,\la\xi,\la y)=\ka_\la\si(x,\xi,y)\ka_\la^{-1},\quad\la\in\RR_+.
\end{equation}
\begin{definition}\label{d18}
The function $\si(x,\xi,y)$ is called the \textit{symbol} of the operator
$A$ with parameter $y\in Y$ on the domain of $\RR^k$ embedded in $X$ and
will be denoted by $\si(A)(x,\xi,y)$.
\end{definition}

\section{Stratified manifolds}

In this section, we describe the class of manifolds on which elliptic
theory will be studied (cf., e.g., \cite{PlSe1,PlSe6,CMSch1}).

\subsection{Definition of stratified manifolds}

Let $\cM$ be a Hausdorff locally compact topological space with a
\textit{filtration of length} $k$, i.e., a decreasing finite chain of
closed subspaces of the form
\begin{equation}\label{filtration}
    \cM\equiv\cM_0\supset\cM_1\supset\dotsm\supset\cM_k.
\end{equation}
The subspaces $\cM_j$ will be called \textit{closed strata}, the sets
\begin{equation*}
    \cM_j^\circ\ovs{\op{def}}=\cM_j\setminus\cM_{j+1}, \quad j=0,\dotsc,k,
\end{equation*}
will be called \textit{open strata},\footnote{For brevity, we adopt the
convention that $\cM_{k+1}=\varnothing$.} and we shall always assume that
the open strata are smooth manifolds without boundary.

The space $\cM$ will be called a \textit{stratified manifold of length} $k$
if the following additional structure is specified on it. The definition of
this structure involves induction over $k$; to define a stratified manifold
of length $k$, one should have the definition not only of a stratified
manifold of length $l\le k-1$ but also of a diffeomorphism of such
manifolds, smoothly depending on parameters, and of the structure of a
stratified manifold of length $l$ on the product of a stratified manifold
of length $l$ by an open subset in $\RR^m$. We start from $k=0$ (the basis
of induction).
\begin{definition}\label{bas-ind}
A \textit{stratified manifold $\cM$ of length zero} is a smooth manifold
without boundary. Diffeomorphisms of smooth manifolds (possibly, smoothly
depending on parameters) and the structure of a smooth manifold on the
product of $\cM$ by an open subset in $\RR^m$ are defined in a standard
way. A \textit{coordinate system of type} $0$ on $\cM$ is a smooth
coordinate system.
\end{definition}
Now we proceed to the general case of a space  $\cM$ with a
filtration~\eqref{filtration}. By analogy with smooth manifolds, we define
a stratified manifold with the use of an \textit{atlas}, i.e., a set of
coordinate systems related by admissible transition maps. On the smooth
part $\cM_0^\circ$ of $\cM$, coordinate systems are the same as in the
smooth case. However, in neighborhoods of points of the singularity set
$\cM_1$ the local models of the manifold are more complicated objects than
open subsets of $\RR^m$. Let $\Om\equiv\Om_0\supset\dotsm\supset\Om_{l}$ be
a stratified manifold of length $l<k$. By $K_\Om$ we denote the infinite
cone
\begin{equation*}
    K_\Om=\{\ov\RR_+\times\Om\}\slash\{\{0\}\times\Om\}
\end{equation*}
with base $\Om$. We denote its vertex by $\mathbf{0}$ and set
\begin{equation*}
 K_\Om^\circ=K_\Om\setminus\{\mathbf{0}\}\equiv\RR_+\times\Om.
\end{equation*}
Local models in a neighborhood of points of $\cM_1$ are open subsets of the
products $\RR^m\times K_\Om$. The ``coordinates'' in such a product will be
denoted by $(x,r,\om)$, where $x\in\RR^m$, $r\in\RR_+$, and $\om\in\Om$.

Let us proceed to precise statements.

A \textit{coordinate system of type} $l$ on $\cM$ (centered at a point
$z\in\cM_l^\circ$) is a smooth coordinate system on $\cM_0^\circ$ in a
neighborhood of $z$ if $l=0$, while for $0<l\le k$ it is a homeomorphism
$\ph:U\lra V$, where
\begin{enumerate}
    \item[1)] $U\subset\cM$ is a neighborhood of $z$ that does not meet
    $\cM_s$ for $l<s\le k$;
    \item[2)] $V$ is neighborhood of the point $(0,\mathbf{0})=\ph(z)$
    in the product $\RR^m\times
    K_\Om$, where $\Om=\Om^{(l)}$ is a compact stratified manifold of
    length $l-1$;
    \item[3)] $\ph(U\cap\cM_l^\circ)=V\cap(\RR^m\times\{\mathbf{0}\})$,
    and the restriction
    $\ph\big|_{U\cap\cM_l^\circ}$, treated as a mapping into
    $\RR^m$, is a smooth coordinate system on $\cM_l^\circ$.
    \item[4)] for $s<l$, one has
    $\ph(U\cap\cM_s^\circ)=V\cap\wt\Om_s$, where
    $\wt\Om_s=\RR^m\times\RR_+\times\Om_s$ are the strata of the stratified
    manifold $\RR^m\times K_\Om^\circ$ of length $l-1$.
\end{enumerate}

Now let $\ph_i:U_i\lra V_i$, $i=1,2$, be two coordinate systems of types
$l_1$ and $l_2$ on $M$ such that $U_1\cap U_2\ne\varnothing$. Then we have
the transition map $\ph_{12}=\ph_1\circ\ph_2^{-1}$. There are two possible
cases: (1) at least one of the numbers $l_i$ is less than $k$; (2)
$l_1=l_2=k$. In the first case, both the domain and the range of $\ph_{12}$
are stratified manifolds of length $l=\min\{l_1,l_2\}<k$ and, using the
induction assumption, we say that the transition mapping $\ph_{12}$ is
\textit{admissible} if it is a diffeomorphism of stratified manifolds of
length $l$. In the second case, we have $\ph_{12}:W\lra\wt W$, where $W,\wt
W\subset\RR^m\times K_\Om$ are some open domains and
$\ph_{12}\br(W\cap(\RR^m\times K_\Om^\circ)\bl)=\wt W\cap(\RR^m\times
K_\Om^\circ)$. We say that $\ph_{12}$ is \textit{admissible} if the
restriction $\ph_{12}\big|_{W\cap(\RR^m\times K_\Om^\circ)}$ is a
diffeomorphism of stratified manifolds of length $k-1$ and in a
sufficiently small neighborhood of the set $\RR^m\times\{\mathbf{0}\}$ the
mapping $\ph_{12}$ has the form
\begin{equation}\label{essence}
    \ph_{12}(x,r,\om)=(\wt\ph(x),r,\psi(x,\om)),
\end{equation}
where $\wt\ph$ is a diffeomorphism and $\psi(x,\bcdot):\Om\lra\Om$ is a
diffeomorphism of a stratified manifolds $\Om$ of length $k-1$ and smoothly
depends on the parameters $x$.

It is easily seen that the composition of admissible transition maps is
again admissible. (The mapping with empty domain is admissible by
definition.)

\begin{definition}\label{step-ind}
The structure of a \textit{stratified manifold of length} $k$ on $\cM$ is a
maximal set of coordinate systems $(U,\ph)$ of types $0,\dotsc,k$ such that
the union of all these $U$ covers $\cM$ and all transition maps between
coordinate systems are admissible.

If $\cM$ is a stratified manifold of length $k$ and $W$ is a domain in
$\RR^m$, then the structure of a stratified manifold $k$ on $W\times\cM$ is
defined as follows: one considers all possible coordinate systems of the
form $\id\times\ph:W\times U\lra W\times V$, where $\ph:U\lra V$ is a
coordinate system on $\cM$, and the resulting atlas is embedded in a
(uniquely determined) maximal atlas.

If $\cM$ and $\cN$ are stratified manifold of length $k$, then a
\textit{diffeomorphism} of $\cM$ onto $\cN$ is a homeomorphism
$f:\cM\lra\cN$ preserving the structure of a stratified manifold (and hence
having the form of admissible transition maps in coordinate systems).

We say that a diffeomorphism $f:\cM\lra\cN$ \textit{smoothly depends on
parameters} if the admissible transition maps representing it in coordinate
systems smoothly depend on these parameters.
\end{definition}
\begin{remark}
Note that admissible transition maps include only diffeomorphisms of
stratified manifolds of length $\le k-1$, for which the notion of smooth
dependence on parameters is known by the induction assumption.
\end{remark}

\begin{remark}
One can readily see that each closed stratum $\cM_j$ of a stratified
manifold of length $k$ is itself a stratified manifold (of length $k-j$).
\end{remark}

\subsection{Generating vector fields of families of diffeomorphisms}

If $f_\tau:\cM\lra\cM$, $\tau\in[0,1]$, $f_0=\id$, is a smooth family of
diffeomorphisms of a stratified manifold $\cM$, then on the dense open set
$\cM^\circ=\cM\setminus\cM_1$ this family defines a vector field $V(z)$,
$z\in\cM^\circ$, (depending on the parameter $\tau$), from which the family
can be reconstructed as the solution of the ordinary differential equation
$dz/d\tau=V(z)$. By our definition of diffeomorphisms, such a vector field
has the following properties:
\begin{enumerate}
    \item[1)] in a coordinate system of type $0$, this is an arbitrary smooth
              vector field;
    \item[2)] in a coordinate system of type $l>0$, for sufficiently small $r$
    the field has the form $V(x,r,\om)=(V_1(x),0,V_2(x,\om))$, where $(V_1(x),V_2(x,\om))$
    is a field with similar properties 1), 2) on the stratified manifold
    $\RR^m\times\Om$ of length $l-1$.
\end{enumerate}
Conversely, each field with these properties generates (at least locally,
and globally on a compact stratified manifold) a family of diffeomorphisms.
The space of vector fields defined by induction on $k$ via properties 1)
and 2) on a stratified manifold $\cM$ of length $k$ will be denoted by
$\cV(\cM)$.

\subsection{The blow-up}

For each stratified manifold, we define its \textit{blow-up}, which is a
continuous mapping
$$
\pi\equiv\pi_\cM: M\longrightarrow \mathcal{M}
$$
of a manifold $M$ with corners onto $\cM$. [Recall that a \textit{manifold
of dimension $n$ with corners} is a Hausdorff space that is locally (i.e.,
in a neighborhood of each point) homeomorphic to the product
$\ov{\mathbb{R}}_+^{k}\times \mathbb{R}^{n-k}$, $0 \le k\le n$, where $k$
can depend on the point, and equipped with an atlas of such homeomorphisms
with transition maps smooth up to the boundary.] We define the blow-up by
the following ``axioms'':
\begin{enumerate}
    \item[1)] If $\cM$ is a stratified manifold of length $0$
    (i.e., a smooth manifold without boundary), then the blow-up has the
    form $\id:\cM\lra\cM$ (i.e., $M=\cM$ and the projection is the identity mapping).
    \item[2)] If $\cN\subset\cM$ is an open (and hence stratified) submanifold, then
    we have the commutative diagram of blow-ups
\begin{equation*}
    \begin{CD}
    N@>>> M\\
    @V{\pi_\cN}VV @VV{\pi_\cM}V\\
    \cN@>>>\cM,
    \end{CD}
\end{equation*}
    where the horizontal arrows are embeddings.
    \item[3)] If $\cM=\RR^m\times K_\Om$, where $\Om$
    is a stratified manifold, then the blow-up has the form
\begin{align*}
    \pi_\cM:\RR^m\times\ov\RR_+\times\wt\Om&\lra\RR^m\times K_\Om,\\
     (x,r,\wt\om)&\longmapsto \begin{cases}
     (x,\mathbf{0}), & r=0\\ \bl(x,r,\pi_\Om(\wt\om)\br),& r>0,
     \end{cases}
\end{align*}
where $\pi_\Om:\wt\Om\lra\Om$ is the blow-up of $\Om$.
\end{enumerate}
One can verify that these axioms by induction give a well-defined blow-up
of every stratified manifold. Moreover, the blow-up is a diffeomorphism of
the interior $M^\circ$ of the manifold $M$ onto the interior
$\cM^\circ\equiv\cM_0^\circ$ of the manifold $\cM$. In what follows, the
term ``blow-up'' will be used for the mapping $\pi:M\lra\cM$ as well as (by
abuse of language) for the manifold $M$ with corners.

\subsection{The cone bundle}

Let $\cM_j^\circ$, $j>0$, be an open stratum of $\cM$. Formally extending
the change of variables~\eqref{essence} written out for coordinate systems
of type $j$ to arbitrarily large $r$, we see that over $\cM_j^\circ$ there is a
well-defined bundle with fiber $K_\Om$, where $\Om\equiv\Om^{(j)}$ is a
compact stratified manifold of length $j-1$.

\begin{remark}
It follows from the construction of the blow-up (we omit the corresponding
computation for lack of space) that this bundle can be canonically (and
smoothly) extended to the \textit{blow-up} $M_j$ of the stratum $\cM_j$.
This remark will be useful in what follows, since the operator-valued
symbols of our $\Psi$DO will be defined over closed strata.
\end{remark}

\subsection{Measures on a stratified manifold}

On a compact stratified manifold $\cM$ of length $k$, we introduce some
natural metrics and the corresponding measures $\mu$, which give rise to
spaces $L^2(\cM)=L^2(\cM,\mu)$ of square integrable functions.

If a Riemannian metric is given, then we have a uniquely determined
Riemannian volume element. Hence it suffices to describe the metrics. This
will be done by induction over the stratification length.

1. On a compact stratified manifold of length 0, we take an arbitrary
Riemannian metric. The corresponding volume element is independent of the
choice of the metric modulo equivalence (i.e., multiplication by an
everywhere positive smooth function).

2. Suppose that we already know how to define metrics and measures on
manifold of length $\le k-1$. To do the same for manifolds of length $k$,
it suffices to consider a neighborhood of the stratum $\cM_k$ (since
outside this neighborhood the stratification length is less than $k$).
Moreover, it suffices to consider neighborhoods in which the cone bundle is
trivialized. For such a neighborhood of the form $U\times K_\Om$, where $U$
is a coordinate neighborhood on $\cM_k$, we define the metric by the formula
\begin{equation}\label{form}
 ds^2=dx^2+dr^2+r^2d\om^2,
\end{equation}
where $d\om^2$ is the metric given by the induction assumption on $\Om$.
Globally, we obtain a metric on $\cM$ by using a partition of unity.

Metrics obtained by this procedure are called \textit{metrics with edge
degeneration}.

Using formula~\eqref{form}, we write out an expression for the measure:
\begin{equation*}
    d\op{vol}=r^ndx\,dr\,d\op{vol}_\Om,
\end{equation*}
where $d\op{vol}_\Om$ is the volume element on $\Om$, known by induction,
and $n=\dim\Om$ is the dimension of $\Om$.

In what follows, all operators on $\cM$ will be considered in the space
$$L^2(\cM)\equiv L^2(\cM,d\op{vol}),$$ and all operators on the cone $K_\Om$
will be considered in the space $$L^2(K_\Om)\equiv L^2(K_\Om,r^n\,
dr\,d\op{vol}_\Om).$$

\subsection{The cotangent bundle}

By induction over the stratification length $k$, we define a space
$\Vect_\cM$ of vector fields on the manifold $\cM$.

1. For a smooth manifold $\cM$, by $\Vect_\cM$ we denote the space of all
vector fields on $\cM$.

2. For a manifold $\cM$ of length $r$, we define $\Vect_\cM$ as the
$C^\infty(M)$-module that is locally (on the set
$U\times\ov\RR_+\times\wt\Om$ corresponding to the coordinate neighborhood
$U\times K_\Om$ of $\cM$) formed by the vector fields
\begin{equation}\label{qqqqqq}
    \theta=a\pd{}v+b\pd{}r+\frac1r\theta_1,
\end{equation}
where $a$ and $b$ are smooth functions and $\theta_1\in\Vect_\Om$.

The metric $ds^2$ defines a $C^\infty(M)$-valued pairing on $\Vect_\cM$.

\begin{proposition}
The formula
\begin{equation*}
    \langle\ph(\theta),\mu\rangle=ds^2(\theta,\mu)
\end{equation*}
specifies a bijection $\varphi$ of the space $\Vect_\cM$ onto a
$C^\infty(M)$-module $\Lambda^1(\cM)\subset\Lambda^1(M)$ of differential
forms on the blow-up $M$ of the manifold $\cM$. The elements of
$\Lambda^1(\cM)$ are exactly the forms than vanish on the fibers of the
projection $p:M\lra\cM$.
\end{proposition}

\begin{proof}
The assertion concerning the bijectivity follows from the fact that
$\Vect_M\subset\Vect_\cM$ and the embedding is epimorphic on the dense main
stratum. The description of elements of $\Lambda^1(\cM)$ follows by a
straightforward computation from~\eqref{qqqqqq}.
\end{proof}

\begin{definition}
The {\em cotangent bundle} $T^*\mathcal{M}$ of the manifold $\cM$ is the
(existing by the Swan theorem) bundle over the blow-up $M$ whose sections
are elements of $\Lambda^1(\cM)$.
\end{definition}

In a similar way, one defines the cotangent bundles of the strata $\cM_j$.

\subsection{The space $C^\infty(\cM)$}

By induction on the length $k$ of a stratified manifold, we define a
function space whose elements are called \textit{smooth functions} on
$\cM$.
\begin{definition}
If $\cM$ is a stratified manifold of length $0$, then the space
$C^\infty(\cM)$ is defined in a standard way.

If $\cM$ is a stratified manifold of length $k>0$, then a function $\ph$ on
$\cM$ belongs to $C^\infty(\cM)$ if and only if
\begin{enumerate}
    \item[1)] $\ph\in C^\infty(\cM\setminus\cM_k)$;
    \item[2)] in each coordinate neighborhood of type $k$ for sufficiently
    small $r$ the function $\ph$ has the form
\begin{equation*}
    \ph(x,r,\om)=\psi(x),
\end{equation*}
where $\psi\in C^\infty(\RR^m)$ is a smooth function.
\end{enumerate}
\end{definition}

It is easily seen that the space $C^\infty(\cM)$ is contained in the space
$C(\cM)$ of continuous functions on $\cM$ and is dense in it. Moreover,
$[\cV(\cM),C^\infty(\cM)]\subset C^\infty(\cM)$.

\section{Pseudodifferential operators}

In this section, we describe the algebra of zero-order $\Psi$DO with
parameters with smooth symbols on a stratified manifold $\cM$. The
definition of $\Psi$DO, as well as the definition of a stratified manifold,
is based on induction over the stratification length. In contrast to the
case of abstract operators considered in Sections 1--4, we shall assume
that the parameter space $Y$ is a finite-dimensional vector space $V$ or,
more generally, a vector bundle with fiber $V$ over a compact smooth
manifold. To avoid cumbersome notation, we usually omit the variables
pertaining to the base of the bundle. It is meant that the estimates
written out below admit as many differentiations as desired along the base.

\subsection{Negligible operators}

The calculus of $\Psi$DO with parameters will be constructed modulo some
space of operator families the addition of which does not affect the
Fredholm property and the index. Here we describe this space. Let $\cM$ be
a stratified manifold.

\begin{definition}
By $J_\infty(V,\cM)\equiv J_\infty(\cM)\equiv J_\infty$ we denote the space
of smooth operator families
\begin{equation}\label{op-sem}
    D(v):L^2(\cM)\lra L^2(\cM),\qquad v\in V,
\end{equation}
such that the following conditions are satisfied:
\begin{enumerate}
    \item[(i)] all operators $D(v)$ are compact in $L^2(\cM)$;
    \item[(ii)] the estimates
\begin{equation}\label{J0}
   \norm{\pd{^\be D(v)}{v^\be}}\le
   C_{\be N}(1+\abs{v})^{-N},\qquad\abs{\be},N=0,1,2,\dotsc\,
\end{equation}
    hold;
    \item[(iii)] conditions (i) and (ii) remain valid if $D(v)$ is replaced by the
    product $$V_1\dotsm V_pD(v)V_{p+1}\dotsm V_{p+q}$$ of arbitrary length $p+q$,
    $p,q\ge0$, where $V_1,\dotsc,V_{p+q}\in\cV(\cM)$.
\end{enumerate}
\end{definition}

\begin{proposition}\label{p5-1}
The space $J_\infty$ is a two-sided $C^\infty(\cM)$-module. If $D\in
J_\infty$ and $f_\a:\cM\lra\cM$ is a family of diffeomorphisms smoothly
depending on parameters $\a$, then the operator family $D_\a\equiv
D_\a(v)=(f_\a^*)^{-1}D(v)f_\a^*$ is an element of $J_\infty$ smoothly
depending on $\a$. \rom(More precisely, the operator $D_\a(v)$ smoothly
depends on $(\a,v)$, and conditions \rom{(i)--(iii)} are satisfied locally
uniformly with respect to $\a$ for the family $D(v)$ itself as well as for
all its derivatives with respect to~$\a$\rom).
\end{proposition}
\begin{proof}
Obviously, conditions (i) and (ii) are preserved under right and left
multiplication by $\ph\in C^\infty(\cM)$. Condition~(iii) is also
preserved, since, as was already noted, $[\cV(\cM),C^\infty(\cM)]\subset
C^\infty(\cM)$, so that it remains to use the Leibniz formula. Finally, let
us prove the smoothness of the family $D_\a$. It depends on the parameter
$\a$ strongly continuously, and the derivative $\pa D_\a/\pa\a_j$ has the
form
\begin{equation*}
    \pd{D_\a}{\a_j}=(f_\a^*)^{-1}[D(v),V]f_\a^*,\quad \text{where $V\in\cV(\cM)$,}
\end{equation*}
on a dense subset in $L^2(\cM)$. By condition~(iii), the element $[D(v),V]$
belongs to $J_\infty$, so that, in particular, $\pa D_\a/\pa\a_j$ lies in
$J_\infty$ and is uniformly bounded and strongly continuous. The same is
true for higher-order derivatives, and hence the function $D_\a$ itself, as
well as its derivatives, is continuous in the operator norm, so that the
desired assertion follows readily.
\end{proof}

\subsection{Definition of pseudodifferential operators}

Let $\cM$ be a compact stratified manifold of length $k$.
\begin{definition}\label{def-PDO}
A \textit{pseudodifferential operator \rom($\Psi$DO) with parameter} $v\in
V$ on $\cM$ is a smooth family of linear operators
\begin{equation}\label{pdo-sem}
    A(v):L^2(\cM)\lra L^2(\cM)
\end{equation}
such that
\begin{enumerate}
    \item[(i)] if $\ph,\psi\in C^\infty(\cM)$ and
    $\supp\ph\cap\supp\psi=\varnothing$, then $\psi A\ph\in J_\infty$;
    \item[(ii)] if $\ph,\psi\in C^\infty(\cM)$ and the supports of both functions
    are contained in some coordinate neighborhood of type $0$, then the operator
    $\psi A\ph$ is a classical zero-order $\Psi$DO in $L^2$ with parameter $v$
    in the sense of Agranovich--Vishik~\cite{AgVi1};
    \item[(iii)] if $\ph,\psi\in C^\infty(\cM)$ and the supports of both
    functions are contained in a coordinate neighborhood of type $l>0$, then the
    operator $\psi A\ph$ is represented in the corresponding coordinate system
    $(x,r,\om)\in\RR^m\times K_\Om$ in the form
\begin{equation}\label{loc-repr}
    \psi A\ph=P\BL(\ovs3x,\ovs3r,\ovs3rv,-i\ovs2{r\pd{}x},
    i\ovs1{r\pd{}r}+i\frac{n+1}2\BR)+R,
\end{equation}
where $R\in J_\infty$, $P(x,r,w,\eta,p)$ is a $\Psi$DO with parameters
$(w,\eta,p)\in\RR^{\dim V+m+1}$ on $\Om$ smoothly depending on the additional
parameters $(x,r)\in\RR^m\times\ov\RR_+$ and compactly supported with
respect to these parameters, and quantization in the first term on the
right-hand side in~\eqref{loc-repr} is understood as follows:
\begin{equation}\label{loc-repr1}
    P\BL(\ovs3x,\ovs3r,\ovs3rv,-i\ovs2{r\pd{}x},
    i\ovs1{r\pd{}r}+i\frac{n+1}2\BR)=\wh P\BL(\ovs2x,-i\ovs1{\pd{}x}\BR),
\end{equation}
where the ``operator-valued symbol'' $\wh P(x,\xi)$ is given by the formula
\begin{equation}\label{loc-repr2}
    \wh P(x,\xi)= P\BL(x,\ovs2r,\ovs2rv,\ovs2{r}\xi,
    i\ovs1{r\pd{}r}+i\frac{n+1}2\BR),
\end{equation}
in which the function of the operator $ir\pa/\pa r$ is defined via the
Mellin transform with the weight line\footnote{This is consistent with the
fact that the $\Psi$DO $P$ is defined for \textit{real} values of the parameter
$p$.} $\im p=-(n+1)/2$.
\end{enumerate}
The set of $\Psi$DO with a parameter $v\in V$ on $\cM$ will be denoted by
$\Psi\equiv\Psi(\cM)\equiv\Psi(V,\cM)$, and by $J=\Psi\cap\cJ$ we denote
its intersection with the ideal $\cJ$ of compact operators decaying as
$v\to\infty$.
\end{definition}

\subsection{The simplest properties of pseudodifferential operators}

In this subsection, we establish some elementary properties of $\Psi$DO. In
particular, we show that they form an algebra in which $J_\infty$ and $J$
are ideals; we also prove some estimates. We have not proved yet that
condition~(iii) in the definition of $\Psi$DO is independent of the
specific choice of a coordinate system (this will be established in the
following subsection), but this does not affect our reasoning.

\begin{theorem}\label{bas-properties}
The following assertions hold.
\begin{enumerate}
    \item[\rom{A.}] A $\Psi$DO $A\in\Psi(V,\cM)$ satisfies the estimates
\begin{equation}\label{esti-psido}
    \norm{\pd{^\a A(v)}{v^\a}}\le
   C_\a(1+\abs{v})^{-|\a|},\qquad\abs{\a}=0,1,2,\dotsc\,,
\end{equation}
    and $A^{(\a)}\in J$ for $\abs{\a}\ge1$.
    \item[\rom{B.}] The set $\Psi(\cM)$ is an algebra, where $C^\infty(\cM)$
    is a subalgebra and $J$ and $J_\infty$ are ideals. Moreover,
\begin{equation*}
    [C^\infty(\cM),\Psi(\cM)]\subset J(\cM).
\end{equation*}
    \item[\rom{C.}] If $V\in\cV(\cM)$, then
    $[V,\Psi(\cM)]\subset\Psi(\cM)$.
\end{enumerate}
\end{theorem}
\begin{proof}
We proceed by induction over the length $k$ of the manifold $\cM$. For
$k=0$, the assertion of the theorem is reduced to standard properties of
$\Psi$DO with a parameter on smooth manifolds, and it suffices to perform
the induction step. It follows from property~(i) that it suffices to verify
everything for operators of the form $\psi A\ph$ localized in some
coordinate system. First, let us prove assertion~A. It follows from
assertion~A known by the induction assumption for operators on manifolds of
length less than $k$, in particular, from the estimates~\eqref{esti-psido},
that the operator family
\begin{equation*}
 F(x,t,v,\xi,p)=P(x,e^{-t},ve^{-t},\xi e^{-t},p)
\end{equation*}
satisfies the estimates
\begin{multline}\label{Opb}
   \norm{\pd{^{\abs{\a}+l+|\be|+|\ga|+k} F(x,t,v,\xi,p)}
   {x^\a\pa t^l\pa v^\be\pa\xi^\ga\pa p^k}}\le
   C_{\a l \be \ga k}(e^t+\abs{v}+\abs{\xi})^{-|\be|-|\ga|}(1+\abs{p})^{-k}\\\le
   \wt C_{\a l \be \ga
   k}(1+\abs{v}+\abs{\xi})^{-|\be|-|\ga|}(1+\abs{p})^{-k},\quad
   \abs{\a}+l+|\be|+|\ga|+k=0,1,2,\dotsc
\end{multline}
and has compact derivatives with respect to the parameters $(v,\xi,p)$. Now
the fact that the operator on the right-hand side in~\eqref{loc-repr} is
well defined, is bounded, and satisfies condition~A is a consequence of the
following general assertion about $\Psi$DO with operator-valued symbols.
\begin{proposition}\label{prdla}
\rom{(i)} Let $H$ be a Hilbert space, and let
$$
H(x,\xi):H\lra H
$$
be a smooth symbol defined on $\RR^{2n}$, ranging in the set of bounded
operators in $H$, and satisfying the estimates
\begin{equation}\label{gensymest}
\norm{\pd{^{\a+\be}H(x,\xi)}{x^\a\pa\xi^\be}}\le
C_{\a\be}(1+\abs{\xi})^{-\abs{\be}},\quad \abs{\a}+
|\be|=0,1,2,\dotsc\,.
\end{equation}
Then the operator
\begin{equation}\label{genoperator}
    \wh H\equiv H\BL(\ovs2x,\smash{-i\ovs1{\pd{}x}}\BR):L^2(\RR^n,H)\lra
    L^2(\RR^n,H)
\end{equation}
is bounded, and its norm is bounded above by a constant that depends only
on finitely many constants $C_{\a\be}$. If, moreover, the symbol is
compact-valued and decays as $(x,\xi)\to\infty$, then the
operator~\eqref{genoperator} is compact.

\rom{(ii)} If the symbol $H(x,\xi)$ is compactly supported with respect to
$x$, then the assertion in~\rom{(i)} remains valid if $H(x,\xi)$ is smooth
in $x$ and only measurable in $\xi$ and the estimates~\eqref{gensymest} are
satisfied only for $\abs{\be}=0$. Under these conditions, if the symbol is
measurable with respect to $\xi$ in the operator norm, is compact-valued,
and uniformly decays as  $\abs{\xi}\to\infty$, then the
operator~\eqref{genoperator} is compact.
\end{proposition}
Indeed, the change of variables $r=e^{-t}$ takes the cone $K_\Om$ to the
cylinder $\Om\times\RR$ and the operator $ir\pa/\pa r$ to $-i\pa/\pa t$. It
remains to note that the operator $-i\pa/\pa t+i(n+1)/2$ is self-adjoint in
the space $L^2$ with weight $e^{-(n+1)t}$ on the cylinder, to which
$L^2(K_\Om)$ is taken by this change of variables, so that the substitution
of this operator as an operator argument for $p$ is meaningful. It remains
to apply Proposition~\ref{prdla} twice.

The proof of B and C is also by induction. It is based on the composition
formula for operators with operator-valued symbol, given in the following
proposition.
\begin{proposition}\label{prdlb}
Let $H_j(x,\xi)$, $j=1,2$, be two symbols satisfying the
estimates~\eqref{gensymest}. Then
\begin{equation}\label{compo}
    \wh H_1\wh H_2=\wh H,
\end{equation}
where the symbol $H(x,\xi)$ also satisfies~\eqref{gensymest} and can be
represented for each $N$ in the form
\begin{equation}\label{product}
    H(x,\xi)=\sum_{\abs{\ga}=0}^{N-1}\frac{(-i)^{|\gamma|}}{|\gamma|!}
    \pd{^\ga H_1(x,\xi)}{\xi^\ga}\pd{^\ga H_2(x,\xi)}{x^\ga}+R_N(x,\xi)
\end{equation}
with a remainder $R_N(x,\xi)$ satisfying the estimates
\begin{equation}\label{ostsymest}
\norm{\pd{^{\a+\be}R_N(x,\xi)}{x^\a\pa\xi^\be}}\le
C_{\a\be}(1+\abs{\xi})^{-N-\abs{\be}},\quad \abs{\a}+
|\be|=0,1,2,\dotsc\,.
\end{equation}
If the $\xi$-derivatives of $H_1$ are compact-valued, then so is
$R_N(x,\xi)$. If, moreover, at least one of the symbols $H_j$ decays as
$x\to\infty$, then the operator $\wh R_N$ is compact.
\end{proposition}
By applying this composition formula to all possible products arising in
the verification of assertions B and C, we arrive at the desired proof.
\end{proof}

\subsection{Invariance with respect to changes of variables}

In this subsection, we show that condition~(iii) in
Definition~\ref{def-PDO} is independent of the choice of a coordinate
system and that $\Psi$DO are invariant under changes of variables.

\begin{theorem}\label{change-var}
The following assertions hold.
\begin{itemize}
    \item[\rom{D.}] If condition~\rom{(iii)} in Definition~\ref{def-PDO}
    is satisfied for the operator $\psi A\ph$ in some coordinate system,
    then it remains valid in any other coordinate system containing
    $\supp\psi\cup\supp\ph$.\footnote{Formulas relating the $\Psi$DO $P$
    in the old and new coordinate systems will be given in the proof.}
    \item[\rom{E.}] If $A\in\Psi(V,\cM)$ and $f_\a:\cM\lra\cM$ is a smooth
    family of diffeomorphisms, then
\begin{equation}\label{change-var-fam}
    A_\a=(f_\a^*)^{-1}A f_\a^*
\end{equation}
is a family of $\Psi$DO with parameters $v\in V$ on $\cM$ smoothly
depending on the additional parameters $\a$.
\end{itemize}
\end{theorem}

\begin{proof}
The assertion of the theorem is known for $k=0$ (smooth manifolds).
Assertion~E follows from D (it suffices to consider a change of coordinates
depending on a parameter), and so we shall prove only~D. It suffices to
consider the case in which both coordinate systems are of type $k$. (If
both systems are of type less than $k$, then the desired result is known by
induction assumption; if only one of the systems is of type $k$, then the
support of the corresponding function $P$ is separated from $r=0$ in this
system, so that the representative in this coordinate system can actually
be assumed to be written out in a coordinate system of type less than $k$).
Next, we can assume that the support of $P$ in both cases is contained in
the domain $r<\e$, where $r$ is sufficiently small, so that the change of
coordinates has the form
\begin{equation*}
    x'=f(x), \om'=g(x,\om).
\end{equation*}
This transformation can be represented as the composition of two
transformations, one of the form $x'=f(x)$ and the other of the form $\om'=
g(x,\om)$. The transformation law for the function $P$ under these changes
of variables can readily be written out. In the first case, we have
\begin{equation}\label{zamenaA}
    P'(x',r,rv,r\xi',p)=P(x,r,rv,r\xi,p),
\end{equation}
where $(x',\xi')=(f(x),(df(x)^*)^{-1}\xi)$, and in the second case
\begin{equation}\label{zamenaB}
    P'(x,0,rv,r\xi,p)=(g_x^*)^{-1}P(x,0,rv,r\xi,p)g_x^*,
\end{equation}
where $g_x=g(x,\bcdot)$. (We omit the more complicated formula for $r>0$,
which should be written out in coordinates on $\Om$, since we do not need
it.) The verification of both formulas is pretty similar to the
finite-dimensional case.
\end{proof}

\section{Symbols and quantization}
The ideal $J$ is contained in the ideal $\cJ$ of compact operators decaying
as the parameter tends to infinity. Hence by  Theorem~\ref{bas-properties},
B $\Psi$DO are local operators with a parameter; i.e.,
$\Psi(\cM)\subset\cA$. Thus we face the problem of computing the local
operators $p_z(A)\in\cA_z$, $z\in\cM$, for a $\Psi$DO $A\in\Psi(\cM)$.

\subsection{Infinitesimal operators and symbols of pseudodifferential operators}

On a stratified manifold, one can define local dilation groups in
coordinates in a neighborhood of each point. Namely, let $z\in\cM_l^\circ$
be some point. Consider a coordinate system of type $l$ centered at $z$.
The point $z$ is depicted by the origin $(0,\mathbf{0})$ in $\RR^m\times
K_\Om$, and we define the action of $\RR_+$ on $\RR^m\times K_\Om$ by the
standard formula
\begin{equation*}
    g_\la(x,r,\om)=(\la x,\la r,\om).
\end{equation*}
The one-parameter unitary group $U_\la$ associated with this action in the
space $L^2(\RR^m\times K_\Om,r^ndxdrd\op{vol}_\Om)=L^2(\RR^m,L^2(K_\Om,r^n
drd\op{vol}_\Om))$, $n=\dim\Om$, is given by the formula
\begin{equation}\label{ulambda}
    [U_\la f](x,r,\om)=\la^{(m+n+1)/2}f(\la x,\la r,\om),
\end{equation}
or
\begin{equation}\label{kalambda}
    [U_\la \mathbf{f}](x)=\la^{m/2}[\ka_\la\mathbf{f}](\la x),
\end{equation}
where the function $\mathbf{f}(x)$ ranges in $L^2(K_\Om,r^n
drd\op{vol}_\Om))$ and $\ka_\la$ is the one-parameter unitary group in the
latter space given by the formula
\begin{equation}\label{kalambda1}
    [\ka_\la u](r,\om)=\la^{(n+1)/2}u(\la r,\om).
\end{equation}
We use this group to construct the infinitesimal operators of $\Psi$DO.
\begin{theorem}\label{PDO-symbols}
The following assertions hold.
\begin{enumerate}
    \item[\rom{F.}] Let $A\in\Psi(V,\mathcal{M})$ be a $\Psi$DO with parameters on $\mathcal{M}$.
    Then for each $z\in\cM_l^\circ$ and each coordinate system of type $l$
    centered at $z$ the operator $A$ has an infinitesimal operator $i_z(A)$
    invariant under $U_\la$. If $\ph,\psi\in C^\infty(\cM)$ are functions
    supported in this coordinate neighborhood and equal to unity at $z$ and
    if $\psi A\ph$ is represented in the local coordinates by
    formula~\eqref{loc-repr}, then the infinitesimal operator is given by
\begin{equation}\label{ifiop}
    i_z(A)=P\BL(0,0,\ovs2rv,-i\ovs2{r\pd{}x},
    i\ovs1{r\pd{}r}+i\frac{n+1}2\BR).
\end{equation}
    \item[\rom{G.}] The infinitesimal operator commutes with translations
    along the stratum $\cM_l^\circ$ (i.e., with respect to the variables $x$),
    and hence the symbol
\begin{equation}\label{ifisym}
    \si(A)(\xi)=P\BL(0,0,\ovs2rv,\ovs2{r}\xi,
    i\ovs1{r\pd{}r}+i\frac{n+1}2\BR)
\end{equation}
    is well defined.
    \item[\rom{H.}] The symbol $\si(A)(\xi)$ smoothly depends on the point
    $x\in\cM_l^\circ$ and is transformed under changes of variables as a function
    on the cotangent bundle of $\cM_l^\circ$ ranging in the space of continuous
    operators in the spaces $L^2$ on the fibers of the cone bundle over $M_l$.
\end{enumerate}
\end{theorem}

\begin{proof}
The proof of assertions F and G is by a straightforward computation with
the use of Propositions \ref{prdla} and \ref{prdlb}. To prove H, we should
additionally use formulas~\eqref{zamenaA} and \eqref{zamenaB}.
\end{proof}

The symbol indicated in the theorem will be denoted by $\si_l(A)(x,\xi)$,
$(x,\xi)\in T^*_0(\cM_l^\circ)$, and called the symbol of $A$ on the $l$th
stratum.

\subsection{Conditions for the Fredholm property of symbols}

By the general localization principle, the Fredholm property (with
parameter) of a $\Psi$DO $A\in\Psi(V,\cM)$ is determined by the
invertibility (for $\xi\ne0$) of its symbols $\si_l(A)(x,\xi)$. Hence
conditions ensuring the invertibility are of interest. Unfortunately, it is
difficult to obtain such conditions in the general case, and so we shall
study a weaker property, namely, the Fredholm property of the symbols.
Consider the symbol~\eqref{ifisym}.
\begin{theorem}\label{fredsym}
The symbol~\eqref{ifisym} is Fredholm for $\xi\ne0$ if and only if
\begin{enumerate}
    \item[\rom{(i)}] all symbols $\si_j(P)(\wt x,\wt \xi)$,
    $j=0,\dotsc,l-1$, of the $\Psi$DO $P(0,0,w,\eta,p)$ with parameters on $\Om$
    are invertible for $\wt \xi\ne0$\rom;
    \item[\rom{(ii)}] the operator family $P(0,0,0,0,p)$ on $\Om$
    is invertible for all $p\in\RR$.
\end{enumerate}
\end{theorem}
\begin{definition}
The symbols $\si_j(P)$ will be called the \textit{symbols} of the symbol
$\si_l(A)$ and will be denoted by $\si_j(\si_l(A))$, $j=0,\dotsc,l-1$. The
family $P(0,0,0,0,p)$ will be called the \textit{conormal symbol} of
$\si_l(A)$ and will be denoted by $\si_c(\si_l(A))$. It is easily seen that
when symbols are multiplied, their symbols (and their conormal symbols) are
also multiplied.
\end{definition}
\begin{proof}[Proof of Theorem~\rom{\ref{fredsym}}]
It follows from Proposition~\ref{prdlb} that if we multiply two
symbols~\eqref{ifisym}, then the corresponding $\Psi$DO $P(0,0,w,\eta,p)$
are also multiplied modulo remainders that give compact-valued symbols.
Hence the assertion of the theorem follows from the fact that
$P(0,0,e^{-t}w,e^{-t}\xi,p)$ is compact-valued and decays as
$\abs{p}\to\infty$ and as $t\to-\infty$ if and only if all symbols of the
family $P$ are zero; it decays also as $t\to\infty$ if and only if the
conormal symbol of $P$ is zero.
\end{proof}

The set of smooth symbols $\si(x,\xi)$ that have the form~\eqref{ifisym}
for each $x$ will be denoted by $\Sigma(T^*_0\cM_l)$.

\subsection{Compatibility conditions and quantization}

So far, we have only established that the symbols of a $\Psi$DO
$A\in\Psi(\cM)$ are defined on the interior of the corresponding cotangent
bundles $T^*_0\cM_j$. In fact, they can be extended by continuity up to the
boundary (where the stratum is adjacent to lower-dimensional strata) and
satisfy certain compatibility conditions on this boundary.

\begin{theorem}\label{soglas}
Let $A\in\Psi(\cM)$ be a $\Psi$DO on $\cM$. Then
\begin{enumerate}
    \item[\rom{K.}] each symbol $\si_j(A)(x,\xi)$ is a smooth function on
    $T^*_0\cM_j$ up to the boundary\rom;
    \item[\rom{L.}] at the points where a stratum $\cM_j$ is adjacent to a stratum $\cM_l$,
    $l>j$, the compatibility conditions
\begin{equation}\label{e-soglas}
    \si_j(D)|_{\cM_l}=\si_j(\si_l(D))
\end{equation}
are satisfied.
\end{enumerate}
\end{theorem}

\begin{proof}
It suffices to compute the symbols of the operator~\eqref{loc-repr} for
$r=r_0>0$. They prove to be equal to the respective symbols of
$P(x,r_0,w,\eta,p)$ and tend to the desired limits as $r_0\to0$.
\end{proof}

By $\Sigma(\cM)\equiv\Sigma(V,\cM)$ we denote the subset of the direct
product
$$
\prod_{j=0}^k\Sigma(T^*_0\cM_j)
$$
formed by elements whose components pairwise satisfy the compatibility
condition~\eqref{e-soglas}. Theorem~\ref{soglas} states that the tuple of
symbols of a pseudodifferential operator $A\in\Psi(\cM)$ always lies in
$\Sigma(\cM)$. It turns out that the converse is also true.

\begin{theorem}\label{antisoglas}
Let $(\si_0,\dotsc,\si_k)\in\Sigma(\cM)$. Then there exists a $\Psi$DO
$A\in\Psi(\cM)$ such that
\begin{equation*}
    \si_j(A)=\si_j,\quad j=0,\dotsc,k.
\end{equation*}
If the symbols $\si_0,\dotsc,\si_k$ smoothly depend on parameters, then $A$
can also be assumed to depend smoothly on the same parameters.
\end{theorem}

\begin{remark}
The operator $A$ is uniquely determined modulo elements of the ideal $J$.
\end{remark}

\begin{proof}
Let us prove this by induction.

1. For operators on a smooth manifold, the existence of an operator with a
given symbol is known.

2. Now suppose that we need to construct an operator with given compatible
symbols on a stratified manifold $\cM$. By subtracting an operator of the
form
$$
\phi(r)\si_0(D)(x,-i\pa/\pa x,v)\phi(r),
$$
where $\phi(r)$ is a cutoff function equal to zero for large $r$, we reduce
the problem to the case in which the symbol corresponding to the stratum
$\cM_k$ is equal to zero identically, and so, by the compatibility
conditions, all other symbols tend to zero when approaching this stratum.

It suffices to construct the corresponding operator on $\cM$ locally, in a
coordinate neighborhood $\RR^m\times K_\Om$ of type $k$. Take a family $\wt
P(x,r,w,\eta,p)$ of $\Psi$DO with parameters $(w,\eta,p)$ on $\Om$
depending on the additional parameters $x,r$ such that $\si_j(P)=\si_j$,
$j=0,\dotsc,k-1$; this is possible by the induction assumption. Since all
these symbols vanish for $r=0$, it follows that the family
\begin{equation*}
    P(x,r,w,\eta,p)=\wt P(x,r,w,\eta,p)-\wt P(x,0,w,\eta,p)
\end{equation*}
has the same symbols $\si_j$ for $j\le k-1$ and vanishes for $r=0$. By
quantizing this family, we obtain the desired operator. The proof of the
theorem is complete.
\end{proof}

\subsection{Main properties of pseudodifferential operators}

\begin{theorem}\label{svoistva-pdo1}
The set $\Sigma({T^*X\times V})$, as well as the set $\Psi(\cM)$ of
$\Psi$DO on a stratified manifold $\cM$, is a local $C^*$-algebra.
\end{theorem}

\begin{proof}
For the norm one can take the supremum of the operator norm over all
parameter values. It has already been proved that these sets are algebras.
Let us show that they are local $C^*$-algebras. To this end, we need to
show that they are closed with respect to holomorphic functional calculus.
Consider, say, the first of these algebras. If an element $D\in
\Sigma({T^*X\times V})$ is invertible, then the inverse also lies in
$\Sigma({T^*X\times V})$. This can readily be derived from the composition
formula. Now let a function $f(z)$ be analytic in a neighborhood of the
spectrum of an element $D$. Then
\begin{equation*}
    f(D)=\frac1{2\pi i}\oint_\ga f(z)(z-D)^{-1}dz,
\end{equation*}
where $\ga$ is a contour surrounding the spectrum and lying in the domain
of analyticity of $f$. By the preceding, the operator $(z-D)^{-1}$ lies in
$\Sigma({T^*X\times V})$. Moreover, one can readily show that this operator
and its symbols continuously depend on the parameter $z$. Hence the
operator$f(D)$ also lies in $\Sigma({T^*X\times V})$, and its symbols are
given by the formula
\begin{equation*}
    \si_j(f(D))=\frac1{2\pi i}\oint_\ga f(z)(z-\si_j(D))^{-1}dz.
\end{equation*}
\end{proof}

Summarizing, we obtain the following theorem.

\begin{theorem}\label{svoistva-pdo2}
The symbol mapping
\begin{equation}\label{symotobr}
    \begin{aligned}
    \si: \Psi(V,\cM) &\lra\bigoplus_{j=0}^k \Sigma({T^*\cM_j\times
    V}) ,\\
    D&\longmapsto (\si_0(D),\dotsc,\si_k(D))
    \end{aligned}
\end{equation}
is a homomorphism of local $C^*$-algebras and generates an isomorphism
\begin{equation*}
    \si:\Psi(V,\cM)\bigm\slash J\lra\Si(V,\cM)\subset \bigoplus_{j=0}^k
    \Sigma({T^*\cM_j\times
    V})
\end{equation*}
onto the local $C^*$-algebra of symbols satisfying the compatibility
conditions~\eqref{e-soglas}.
\end{theorem}

\subsection{Ellipticity and the Fredholm property}

\begin{definition}
An operator $D\in\Psi(V,\cM)$ is said to be \textit{elliptic} if all
symbols $\si_j(D)$, $j=0,\dotsc,k$, are invertible outside the zero
sections of the corresponding bundles.
\end{definition}

The general localization principle implies the finiteness theorem in a
standard way.
\begin{theorem}\label{ellipticity}
Elliptic $\Psi$DO on a compact stratified manifold $\cM$ are Fredholm with
a parameter. In particular, if $V\ne\{0\}$, then an operator elliptic with
a parameter is invertible for large $\abs{v}$.
\end{theorem}

\section{Proof of Proposition~\ref{auxp1}}

Consider the composition
\begin{equation*}
    B=\cF_{t\to\xi}\wh B\ov\cF_{\xi\to t}.
\end{equation*}
We wish to show that $B$ can be represented as the pointwise action of an
operator-valued function $B(\xi)$ with the properties indicated in the
proposition. Since the operator $\wh B$ commutes with translations, it
follows that $B$ commutes with the multiplications by the exponentials
$\exp (it\xi)$, $t\in\RR^k$, and hence, by continuity, with the
multiplication by an arbitrary function $\ph(\xi)\in C_0^\infty(\RR^k)$.
Now let $e\in H$ be an arbitrary vector, and let $\chi(\xi)\in
C_0^\infty(\RR^k)$ be an arbitrary function equal to unity on the support
of $\ph$. Then
\begin{equation}\label{zagadka}
    B (\ph e)=B(\ph\chi e)=\ph B(\chi e).
\end{equation}
The function $[B(\chi e)](\xi)$ belongs to $L^2(\RR^k,H)$ and, in
particular, is measurable. Relation~\eqref{zagadka} can be treated as
equality a.e.\ of two measurable $H$-valued functions of $\xi$. At the
points $\xi$ such that $\ph(\xi)\ne0$, one has a.e.
\begin{equation*}
 [B(\chi e)](\xi)=\ph^{-1}(\xi)B(\ph e)(\xi),
\end{equation*}
where the right-hand side is independent of $\chi$ and the left-hand side
is independent of $\ph$. Since the compactly supported functions $\ph(\xi)$
and $\chi(\xi)$ are subjected to the only condition $\ph\chi=\ph$, we
readily see that for each function $\ph(\xi)\in C_0^\infty(\RR^k)$ one has
a.e.
\begin{equation*}
    [B (\ph e)](\xi) = \ph(\xi) f(\xi),
\end{equation*}
where $f(\xi)$ is a measurable $H$-valued function. A simple argument shows
that
\begin{equation*}
    \esssup\norm{f(\xi)}\le\norm{B}=\norm{\wh B}.
\end{equation*}
If $e=e_1+e_2$, then the corresponding functions $f(\xi)$ a.e.\ satisfy
$f(\xi)=f_1(\xi)+f_2(\xi)$ (since this holds, by the linearity of the
operator, for the functions $\ph f$). Now let $\gZ\in H$ be a dense
additive subgroup (say, the set of finite linear combinations of vectors of
some basis with rational coefficients). For each $e\in\gZ$, we fix a
representative of the corresponding element $f\in L^2(\RR^k,H)$. Let us
introduce the sets
\begin{equation*}
    \Delta_e=\{\xi\mid\norm{f(\xi)}\ge\norm{B}\},\quad\Delta_{e_1e_2}
    =\{\xi\mid f(\xi)\ne f_1(\xi)+f_2(\xi),\text{ where }e=e_1+e_2\}.
\end{equation*}
They are of measure zero, and so is their union
\begin{equation*}
    \Delta=\bigcup_{e\in\gZ}\Delta_e\cup\bigcup_{e_1,e_2\in\gZ}\Delta_{e_1e_2}.
\end{equation*}
For each $\xi\in\RR^k\setminus\Delta$, the mapping $e\mapsto f(\xi)$
defined on $\gZ$ is linear and is bounded by the norm $\norm{B}$; hence its
closure is a bounded linear operator $B(\xi):H\lra H$. Set $B(\xi)=0$ for
$\xi\in\Delta$. The operator function this defined is measurable (since it
is uniformly bounded and measurable on the dense set $\gZ$) and satisfies
all desired conditions.

Finally, the continuous dependence of the operator $B(\xi)$ on the
additional parameters is obvious. Indeed, by the preceding we have
\begin{equation*}
    \norm{B_\tau(\xi)-B_{\tau_0}(\xi)}\le\norm{\wh B_y-\wh B_{y_0}}\to0.\hfill\qed
\end{equation*}

\end{document}